\documentclass[11pt]{article}

\usepackage{amsmath,amssymb,amsfonts,enumitem,float,xy}
\xyoption{all}
\usepackage{pstricks,pst-node}

\numberwithin{equation}{section}
\newtheorem{thm}{Theorem}[section]
 
\newtheorem{lmm}[thm]{Lemma}  
\newtheorem{cnj}{Conjecture}  
\newtheorem{cnd}{Condition}

\newcounter{pr}

\newcounter{zr}

\def\BE#1{\begin{equation}\label{#1}}
\def\EE{\end{equation}}
\def\lan{\langle} \def\ran{\rangle}

\def\blr#1{\big\lan{#1}\big\ran}
\def\llrr#1{\lan\!\lan{#1}\ran\!\ran}
\def\bllrr#1{\big\lan\!\big\lan{#1}\big\ran\!\big\ran}
\def\e_ref#1{(\ref{#1})}
\def\ov#1{\overline{#1}} \def\un#1{\underline{#1}}
\def\ti#1{\tilde{#1}} \def\wt#1{\widetilde{#1}}

\def\sm#1{\begin{small}#1\end{small}}
\def\Lg#1{\begin{Large}#1\end{Large}}
\def\Hg#1{\begin{Huge}#1\end{Huge}}

\def\QED{$\Box$}
\def\lra{\longrightarrow}

\def\be{\beta}
\def\de{\delta}
\def\ep{\epsilon}
\def\ga{\gamma}
\def\io{\iota}
\def\ka{\kappa}
\def\la{\lambda}
\def\om{\omega}
\def\si{\sigma}

\def\ups{\upsilon}
\def\ve{\varepsilon}
\def\vph{\varphi}

\def\De{\Delta}
\def\Ga{\Gamma}
\def\Si{\Sigma}

\def\C{\mathbb C}
\def\cC{\mathcal C}
\def\tC{\tt C}

\def\tC{\mathtt C}
\def\cD{\mathcal D}
\def\D{\mathfrak D}
\def\d{\mathfrak d}
\def\E{\mathbb E}

\def\cM{\mathcal M}
\def\M{\mathfrak M}
\def\N{\mathcal N}
\def\cO{\mathcal O}
\def\fO{\mathfrak O}
\def\P{\mathbb P^1}
\def\bP{\mathbb P}
\def\Q{\mathbb Q}
\def\R{\mathbb R}
\def\cS{\mathcal S}
\def\T{\mathcal T}
\def\U{\mathfrak U}
\def\Z{\mathbb Z}
\def\cZ{\mathcal Z}

\def\i{\infty}
\def\dbar{\bar\partial}
\def\eset{\emptyset}

\def\Aut{\textnormal{Aut}}

\def\Cont{\textnormal{Cont}}
\def\dim{\textnormal{dim}}

\def\Edg{\textnormal{Edg}}
\def\eff{\textnormal{eff}}
\def\ev{\textnormal{ev}}
\def\gh{\textnormal{gh}}
\def\Hom{\textnormal{Hom}}
\def\Im{\textnormal{Im}}
\def\mod{\textnormal{mod}\,}
\def\PD{\textnormal{PD}}
\def\rk{\textnormal{rk}}

\def\Ver{\textnormal{Ver}}

\begin{document}

\title{Enumerative Geometry of Calabi-Yau $5$-Folds}
\author{R. Pandharipande and A. Zinger}
\date{11 February 2008}

\maketitle

\begin{abstract}

Gromov-Witten theory is used to define an enumerative geometry of 
curves in Calabi-Yau $5$-folds.
We find 
recursions for meeting numbers of genus 0 curves, and
we determine the contributions of moving multiple 
covers of genus 0 curves to the genus $1$ Gromov-Witten invariants. 
The resulting invariants, conjectured to
be integral, are analogous to the previously defined
BPS counts for Calabi-Yau $3$ and $4$-folds.
We comment on the situation in higher dimensions
where new issues arise.

Two main examples are considered: the local Calabi-Yau  
$\bP^2$ with normal bundle $\oplus_{i=1}^3\cO(-1)$ and
the compact Calabi-Yau  hypersurface $X_7 \subset\bP^6$.
In the former case, a closed form for our integer
invariants has been conjectured by G. Martin. In the
latter case, we recover in low degrees the classical
enumeration of
elliptic curves by Ellingsrud and Str\"omme.

\end{abstract}

\tableofcontents

\setcounter{section}{-1}
\section{Introduction}
\label{intro_sec}

\subsection{Overview}
\label{overview_subs}

\noindent
Let $X$ be a nonsingular projective variety over $\C$.
Let $\ov\M_{g,k}(X,\be)$ be the moduli space of genus $g$,
$k$ pointed stable maps to $X$ representing the
class $\be\in H_2(X,\Z)$.
Let
$$\ev_i\!: \ov\M_{g,k}(X,\be)\lra X$$
be the evaluation morphism at the $i^{th}$ marking.
The Gromov-Witten theory of primary fields concerns the invariants
\BE{GWdfn_e}
N_{g,\be}(\gamma_1,\ldots,\gamma_k)
= \int_{[\ov\M_{g,k}(X,\be)]^{vir}}
\prod_{i=1}^{k}\ev_i^*(\gamma_i) \ \in \mathbb{Q},\EE
where $\gamma_i\!\in\!H^*(X,\Z)$.
The relationship between the Gromov-Witten invariants 
and the actual enumerative geometry of curves in $X$
is subtle. An overview
of the subject in low dimensions can be found in the introduction of
 \cite{KlP}.

For Calabi-Yau $3$-folds, the Aspinwall-Morrison formula \cite{AsM}
is conjectured to produce integer invariants in genus $0$.
A full integrality conjecture for the Gromov-Witten theory of 
Calabi-Yau $3$-folds was formulated by Gopakumar and Vafa in \cite{GV1,GV2}
in terms of BPS states with geometric
 motivation partially provided by \cite{degencontr}.
The Aspinwall-Morrison prediction has been
extended to all Calabi-Yau $n$-folds in \cite{KlP}:
the numbers $n_{0,\be}(\gamma_1,\ldots,\gamma_k)$ defined by
\BE{g0BPS_e}
\sum_{\be\ne0}N_{0,\be}(\gamma_1,\ldots,\gamma_k)q^{\be}
=\sum_{\be\ne0}n_{0,\be}(\gamma_1,\ldots,\gamma_k)
\sum_{d=1}^{\i}\frac{1}{d^{3-k}}q^{d\be}\EE
are conjectured to be integers.

Let $X$ be a Calabi-Yau of dimension $n\geq 4$.
Since Gromov-Witten invariants of genus $g\!\ge\!2$ of $X$
vanish for dimensional reasons,
only integrality predictions for 
genus $1$ invariants of $X$
remain to be considered.
The analogue of the genus 1 Gopakumar-Vafa integrality prediction
for Calabi-Yau $4$-folds has been formulated in
\cite{KlP}.
Here, we find complete formulas in 
dimension $5$ and reinterpret the dimension 4
predictions.  The geometry becomes significantly more complicated
in each dimension. We discuss new aspects
of the higher dimensional cases.

The relationship between Gromov-Witten theory and enumerative
geometry in dimensions greater than 3 is simplest in the
Calabi-Yau case. The Fano case, even in dimension 4,  
involves complicated higher genus phenomena which 
have not yet been understood.

\subsection{Elliptic invariants}
If $X$ is Calabi-Yau, the virtual moduli cycle for 
$\ov\M_1(X,\be)$ is of dimension $0$.
We denote the associated Gromov-Witten invariant by $N_{1,\be}$,
$$N_{1,\be}=\int_{[\ov\M_1(X,\be)]^{vir}}1 ~\in\Q.$$

Integrality predictions for Calabi-Yau $n$-folds are obtained
by relating curve counts to Gromov-Witten invariants in 
an ideal Calabi-Yau $X$.
All  genus $1$ curves in $X$ are assumed to be  nonsingular, 
super-rigid{\footnote{A nonsingular curve $E\subset X$ with normal
bundle $\N_E$ is super-rigid if,
for every dominant stable map $f:C\rightarrow E$, the vanishing  
$H^0(C,f^*\N_E)=0$ holds.},
and  disjoint from other curves.
Each genus $1$ degree $\be$ curve then contributes $\si(d)/d$ 
to $N_{1,d\be}$ for every $d\!\in\!\Z^+$
via \'etale covers, where 
$$\si(d)=\sum_{i|d} i.$$
The genus 1 to genus 1 multiple cover contribution is independent
of dimension.

If $X$ is an ideal Calabi-Yau $3$-fold, 
the genus $0$ curves in $X$ are also nonsingular, super-rigid, and
disjoint.
The contribution of a genus 0 degree $\be$ curve to $N_{1,d\be}$
is then the integral of an Euler class of an
obstruction bundle on $\ov\M_1(\P,d)$,
$$\int_{[\ov\M_1(\P,d)]^{vir}} e(\text{Obs}) = \frac{1}{12d},$$ 
calculated in \cite{degencontr}.
Thus, if $X$ is an ideal Calabi-Yau $3$-fold,
\BE{g1n3dfn_e}
\sum_{\be\ne0}N_{1,\be}q^{\be}
=\sum_{\be\ne0}n_{1,\be}\sum_{d=1}^{\i}\frac{\si(d)}{d}q^{d\be}
-\frac{1}{12}\sum_{\be\ne0}n_{0,\be}\log(1-q^{\be})\,,\EE
where the enumerative invariant $n_{1,\be}$ is 
defined by \eqref{g1n3dfn_e} and the genus 0 invariant
$n_{0,\be}$ is defined by the 
Aspinwall-Morrison formula \eqref{g0BPS_e}.
The invariants $n_{1,\beta}$ are then conjectured to be integers for
all Calabi-Yau $3$-folds.

If $X$ is an ideal Calabi-Yau $4$-fold, embedded
genus 0 degree $\be$ curves in $X$
form a nonsingular, compact, $1$-dimensional family $\ov\cM_{\be}$.
The moving multiple cover calculation of
 Section 2 of \cite{KlP} 
shows that
$\ov\cM_{\be}$ contributes 
$\chi(\ov\cM_{\be})/24d$ to $N_{1,d\be}$ for every $d\!\in\!\Z^+$.
The calculation is done in two steps. First, the
moving multiple cover integral is done
 assuming every genus 0 degree $\be$ curve is nonsingular.
Second, the contribution from the nodal curves is determined
for a particular, but sufficiently representative, Calabi-Yau $4$-fold $X$
by localization.
For an ideal Calabi-Yau $4$-fold $X$,
\BE{g1n4dfn_e}
\sum_{\be\ne0}N_{1,\be}q^{\be}
=\sum_{\be\ne0}n_{1,\be}\sum_{d=1}^{\i}\frac{\si(d)}{d}q^{d\be}
-\frac{1}{24}\sum_{\be\ne0}\chi(\ov\cM_{\be})\log(1-q^{\be})\,.\EE
The topological Euler characteristic $\chi(\ov\cM_{\be})$  is
determined by
$$\chi(\ov\cM_{\be}) = -n_{0,\be}(c_2(X)) + \sum_{\be_1+\be_2=\be} 
m_{\be_1,\be_2},$$
where $m_{\be_1,\be_2}$ is
the number of ordered pairs $(\cC_1,\cC_2)$  of rational
curves of classes~$\be_1$ and~$\be_2$
meeting at point, see Section 1.2 of \cite{KlP}.

The meeting numbers $m_{\be_1,\be_2}$ can be expressed 
in terms of the invariants $n_{0,\be}(\gamma)$ through a recursion 
on the total degree $\be_1\!+\!\be_2$ by computing the excess contribution
to the topological Kunneth decomposition of $m_{\be_1\be_2}$,
see Sections~0.3 and~1.2 of~\cite{KlP}.
Along with these recursions, relations \e_ref{g0BPS_e} 
and \e_ref{g1n4dfn_e} effectively determine 
the numbers $n_{1,\be}$ in terms of the genus~0 and genus~1 
Gromov-Witten invariants of $X$. 
For arbitrary Calabi-Yau $4$-folds, equation \eqref{g1n4dfn_e}
is taken to be the definition of the numbers $n_{1,\be}$
which are conjectured always to be integers.

If $X$ is an ideal Calabi-Yau $5$-fold,  embedded
genus 0 degree $\be$ curves in $X$ 
form a nonsingular, compact, $2$-dimensional family~$\ov\cM_{\be}$.
However, as the nodal curves are more complicated, the localization
strategy of \cite{KlP} does not appear possible.
By viewing $N_{1,d\be}$ as the number of solutions, 
counted with appropriate multiplicities, of 
a perturbed $\dbar$-equation as in \cite{FO,LT},
we show in Section~\ref{g1nums_sec} that $\ov\cM_{\be}$ contributes 
$$\frac{1}{24d}\,\int_{\ov\cM_{\be}}
\big(2c_2(\ov\cM_{\be})\!-\!c_1^2(\ov\cM_{\be})\big)$$
to $N_{1,d\be}$ for every $d\!\in\!\Z^+$.
Thus, for an ideal Calabi-Yau $5$-fold~$X$,
\BE{g1gendfn_e}\begin{split}
\sum_{\be\ne0}N_{1,\be}q^{\be}
=\sum_{\be\ne0}n_{1,\be}\sum_{d=1}^{\i}\frac{\si(d)}{d}q^{d\be}
-\frac{1}{24}\sum_{\be\ne0}
\int_{\ov\cM_{\be}}\big(2c_2(\ov\cM_{\be})\!-\!c_1^2(\ov\cM_{\be})\big)
 \cdot \log(1-q^{\be})\,.\end{split}\EE

The last term in \e_ref{g1gendfn_e} may be  written 
in terms of various 
meeting numbers of total degree $\be$ 
via a Grothendieck-Riemann-Roch computation
applied to the deformation characterization of the tangent
bundle $T\ov\cM_{\be}$.
We pursue a more efficient strategy in Sections~\ref{g0nums_sec} and~\ref{g1nums_sec}.
Degree $1$ maps from genus 0 curves to degree $\be$ curves in $X$ 
are regular.
Thus, equation~(2.15) 
in \cite{g1diff} expresses their contribution to $N_{1,\be}$ in terms of
counts of $m$-tuples of 1-marked curves with cotangent $\psi$-classes 
meeting at the marked point. 
The $\psi$-classes can be easily eliminated using 
the topological recursion relation at the cost of introducing 
counts of arbitrary meeting configurations of rational curves in~$X$.
The latter can be recursively defined as in the case
of $m_{\be1,\be2}$ in dimension~4.
Relations \e_ref{g0BPS_e} and \e_ref{g1gendfn_e} 
then reduce the numbers $n_{1,\be}$ to functions of 
genus 0 and genus 1 Gromov-Witten invariants. 

Let $X$ be an arbitrary Calabi-Yau $5$-fold. Equation \eqref{g1gendfn_e}
together with the rules provided in Sections
\ref{g0nums_sec} and \ref{g1nums_sec} for the calculation of
$$\int_{\ov\cM_{\be}}\big(2c_2(\ov\cM_{\be})\!-\!c_1^2(\ov\cM_{\be})\big)$$
in terms of the Gromov-Witten invariants of $X$
{\em define} the invariants $n_{1,\be}$. We view 
$n_{1,\beta}$ as virtually enumerating elliptic curves in $X$.

\begin{cnj} For all Calabi-Yau 5-folds $X$ and 
curve classes $\beta\neq 0$, the invariants
$n_{1,\beta}$ are integers.
\end{cnj}


\subsection{Examples}
If the Gromov-Witten invariants of $X$ are known, 
equation \eqref{g1gendfn_e} provides an effective
determination of the elliptic invariants $n_{1,\be}$.
We consider two representative examples.

The most basic local Calabi-Yau 5-fold is the total
space of the bundle
\begin{equation}
\label{xxcc}
\cO(-1)\oplus\cO(-1)\oplus\cO(-1)\lra \bP^2.
\end{equation}
The balanced property of the bundle is
analogous to the fundamental local Calabi-Yau
3-fold
$$\cO(-1)\oplus\cO(-1)\lra \bP^1.$$
As in the 3-fold case, we find very simple closed forms in
Section \ref{lll1}
for the genus 0 and 1 Gromov-Witten invariants of 
the local Calabi-Yau 5-fold \eqref{xxcc}.

We have computed the invariants $n_{1,d}$ via equation \eqref{g1gendfn_e}
up to  degree~$200$. All are integers.
Even the first $60$, shown in Table \ref{g1local_table},
suggest intriguing patterns. For example, 
$n_{1,d}\!=\!0$ for all multiples of $8$.
G. Martin has proposed an explicit formula for $n_{1,d}$
which holds for all the numbers we have computed.
We state Martin's conjecture in Section \ref{lll2}.

\begin{table}
\begin{center}
\begin{tabular}{|r|r||r|r||r|r||r|r||r|r||r|r|}
\hline
$d$& $n_{1,d}$& $d$& $n_{1,d}$& $d$& $n_{1,d}$& $d$& $n_{1,d}$& $d$& $n_{1,d}$&
$d$& $n_{1,d}$\\
\hline
1&    0&    11&   -225&    21&   3025&  31& -14400&  41& -44100& 51&  105625\\
2&    0&    12&    -19&    22&   3870&  32&      0&  42& -51590& 52&   -7119\\
3&   -1&    13&   -441&    23&  -4356&  33&  18496&  43& -53361& 53& -123201\\
4&    0&    14&    630&    24&      0&  34&  22140&  44&  -3645& 54&       0\\
5&   -9&    15&    784&    25&      0&  35&  23409&  45&      0& 55&  142884\\
6&   20&    16&      0&    26&   7560&  36&      0&  46&  74250& 56&       0\\
7&  -36&    17&  -1296&    27&      0&  37& -29241&  47& -76176& 57&  164836\\
8&    0&    18&      0&    28&   -594&  38&  34560&  48&      0& 58&  187740\\
9&    0&    19&  -2025&    29& -11025&  39&  36100&  49&      0& 59& -189225\\
10& 162&    20&   -153&    30& -13412&  40&      0&  50&      0& 60&   12628\\
\hline
\end{tabular}
\end{center}
\caption{Invariants $n_{1,d}$ for $\cO(-1)\oplus\cO(-1)\oplus\cO(-1)\lra\bP^2$}
\label{g1local_table}
\end{table}

The Calabi-Yau septic hypersurface $X_7 \subset \bP^6$
is a much more complicated example.
Using the closed formulas for the genus $1$ and $2$-pointed genus $0$
Gromov-Witten
invariants provided by \cite{g1diff} and \cite{bcov0} respectively,
we have computed  $n_{1,d}$ for $d\!\le\!100$.
All are integers.
The values of $n_{1,d}$ for $d\leq 10$  are shown in Table~\ref{g1_table}.

The invariants $n_{1,d}$ for $d\leq 4$ agree with known enumerative
results for $X_7$. The invariants $n_{1,1}$ and $n_{1,2}$ vanish
by geometric considerations.
Since every genus $1$ curve of degree $3$ in $\bP^6$ is planar,
the number of elliptic cubics on a general $X_7$ can be computed 
classically via Schubert calculus.
The classical calculation agrees with $n_{1,3}$.
Using the expression of non-planar genus 1 curves of degree~4 as
complete intersections of quadrics, Ellingsrud and Str\"omme
have enumerated elliptic quartics on $X_7$ in
Theorem 1.3 of \cite{ES}.
The result agrees with $n_{1,4}$.
To our knowledge, the numbers $n_{1,d}$ are inaccessible
by classical techniques for $d\geq 5$.

\begin{table}
\begin{center}
\begin{tabular}{|r|r|}
\hline
$d$& $n_{1,d}$\\
\hline
1& 0\\
2& 0\\
3& 26123172457235\\
4& 81545482364153841075\\
5& 117498479295762788677099464\\
6& 126043741686161819224278666855602\\
7& 117293462422824431122974865933687206294\\
8& 100945295955344375879041227482174735213546636\\
9& 82898589348613625712387472944689576403215969839772\\
10& 66074146583335641807745540088333857250772567526848951526\\
\hline
\end{tabular}
\end{center}
\caption{Invariants $n_{1,d}$ for a degree $7$ hypersurface in $\bP^6$}
\label{g1_table}
\end{table}

\subsection{BPS states}

The 
integer expansion \eqref{g1gendfn_e}
can be alternatively written
as
\BE{g1gendfn_e2}\begin{split}
\sum_{\be\ne0}N_{1,\be}q^{\be}
=-\sum_{\be\ne0}\tilde n_{1,\beta}\cdot \log(1-q^\beta)
-\frac{1}{24}\sum_{\be\ne0}
\int_{\ov\cM_{\be}}\big(2c_2(\ov\cM_{\be})\!-\!c_1^2(\ov\cM_{\be})\big)
 \cdot \log(1-q^{\be})\,.\end{split}\EE
The integrality condition for the invariants
$\tilde{n}_{1,\beta}$ is equivalent to the conjectured
integrality for ${n}_{1,\beta}$.
We 
view the invariants $\tilde n_{1,\beta}$ as analogous to the BPS state
counts in dimensions 3 and~4.

\subsection{Higher dimensions}
The family $\ov\cM_{\be}$ of embedded genus~$0$ degree~$\be$ 
curves in~$X$ is nonsingular and compact for
 ideal Calabi-Yau $n$-folds for
 $n\!=\!3,4,5$. The moving multiple cover results for
$n\!=3,4,5$ can be summarized by the following equation. 
The contribution of $\ov\cM_{\be}$
to the genus $1$ degree $d\be$ Gromov-Witten invariant
is 
\BE{scaling_e}
\tC_{\be}(d\be)=\frac{1}{24d}\, \int_{\ov\cM_{\be}}
\big(2c_{n-3}(\ov\cM_{\be})\!-\!c_1(\ov\cM_{\be})c_{n-4}(\ov\cM_{\be})\big)
\,.\EE
For dimension 6 and higher, the family of embedded
genus $0$ degree $\be$ curves in $X$
is not compact (multiple covers can occur as limits) even
in ideal cases. 
Nevertheless, we expect a  contribution equation of
the form of \eqref{scaling_e} to hold.
The result
should yield integrality predictions in higher dimensions.

Since the complexity of the Gromov-Witten approach 
increases so much in every dimension, an alternate method
for dimensions 6 and higher is preferable. It is hoped
a connection to newer sheaf enumeration  and derived category
techniques will be made \cite{PT,RT}.

\subsection{Acknowledgments}
We thank J. Bryan, I. Coskun, A. Klemm, J. Starr, and
R. Thomas for several related
discussions. 
We are grateful to G. Martin for finding the pattern governing
the invariants $n_{1,d}$ for the local Calabi-Yau 5-fold geometry. 

The research was started during a visit to the Centre de Recherches 
Math\'ematique in Montr\'eal in the summer of 2007.
R.P. was partial supported by DMS-0500187. 
A.Z. was partially supported by the Sloan foundation and DMS-0604874.

\section{Genus 0 invariants}
\label{g0nums_sec}

\subsection{Configuration spaces of genus 0 curves}
\label{g0config_subs}

\noindent
Let $X$ be a Calabi-Yau 5-fold.
We specify here what conditions an ideal $X$ 
is to satisfy with respect to genus  0 curves.
We denote by
$$H_+(X)\subset H_2(X,\Z)-0$$
the cone of effective curve classes.
If $\be,\be'\!\in\!H_+(X)$, we write $\be'\!<\!\be$ if 
$\be\!-\!\be'$ is an element of~$H_+(X)$.

If $J$ is a finite set and $\be\!\in\!H_+(X)$, we denote by 
$\ov\M_{0,J}(X,\be)$ the moduli space of genus $0$, $J$-marked
stable maps to $X$ representing the class $\be$. 
For $j\!\in\!J$, let 
$$L_j\lra\ov\M_{0,J}(X,\be)$$
be the universal tangent line bundle at the $j$th marked point.
Denote by
$$\cD_j\in\Ga\big(\ov\M_{0,J}(X,\be),\Hom(L_j,\ev_j^*TX)\big)$$
the bundle section induced by the differential of the stable maps
at the $j^{th}$ marked point.

If $\Si$ is a curve, a map $u\!:\Si\!\lra\!X$ is called {\em simple} 
if $u$ is injective on the complement of finitely many points and of
the components of~$\Si$ on which $u$ is constant.
We will call a tuple $(u_1,\ldots,u_m)$ of maps $u_i\!:\Si\!\lra\!X$
{\em simple} if the map
$$\bigsqcup_{i=1}^{m}\Si_i\lra X, \qquad
z\lra u_i(z)~~\hbox{if}~~z\in\Si_i,$$
is simple.
If $J$ is a finite set and $\be\!\in\!H_+(X)$, let
$$\M_{0,J}^*(X,\be)\subset \ov\M_{0,J}(X,\be)$$
be the open subspace of stable maps $[\Si,u]$ such that 
$\Si$ is a $\bP^1$ and $u$ is a simple map.

If $J_1$ and $J_2$ are two finite sets and $\be_1,\be_2\!\in\!H_+(X)$, 
we denote by
$$\M_{0,(J_1,J_2)}^*\big(X,(\be_1,\be_2)\big)
\subset\big\{(b_1,b_2)\!\in\!
\M_{0,\{0\}\sqcup J_1}^*(X,\be_1)\!\times\!\M_{\{0\},0\sqcup J_2}^*(X,\be_2)\!:
\ev_0(b_1)\!=\!\ev_0(b_2)\big\}$$
the subset of simple pairs of maps.
Similarly, if $\be_1,\be_2,\be_3\!\in\!H_+(X)$, let
$$\M_{0,\eset}^*\big(X,(\be_1,\be_2,\be_3)\big)
\subset \big\{(b_1,b_2,b_3)\!\in\!
\M_{0,(\eset,\{1\})}^*(X,(\be_1,\be_2))\!\times\!\M_{0,\{0\}}^*(X,\be_3)\!:
\ev_1(b_2)\!=\!\ev_0(b_3)\big\}$$
be the subset of simple triples of maps.
If $X$ is an ideal Calabi-Yau $5$-fold satisfying
Conditions \ref{g0reg_cnd} and \ref{g0ev_cnd} below, 
there are no other configurations of simple genus 0 curves in $X$,
see Figure~\ref{dim5conf_fig}.

Denote by $\ov\M_{0,J}^*(X,\be)\subset \ov\M_{0,J}(X,\be)$ and
$$\ov\M_{0,(J_1,J_2)}^*\big(X,(\be_1,\be_2)\big)\subset
\ov\M_{0,\{0\}\sqcup J_1}(X,\be_1)\!\times\!\ov\M_{\{0\},0\sqcup J_2}(X,\be_2),$$
the closures of $\M_{0,J}^*(X,\be)$ and 
$\M_{0,(J_1,J_2)}^*\big(X,(\be_1,\be_2)\big)$.
Let
\BE{pi12_e}
\pi_1,\pi_2\!:\ov\M_{0,(J_1,J_2)}^*\big(X,(\be_1,\be_2)\big)\lra
\ov\M_{0,\{0\}\sqcup J_1}\big(X,\be_1\big), 
\ov\M_{0,\{0\}\sqcup J_2}\big(X,\be_2\big),\EE
be the component projection maps.

\begin{figure}
\begin{pspicture}(-3.5,-2.2)(10,1.25)
\psset{unit=.4cm}
\psarc(0,-1){3}{120}{240}\rput(-1.2,1){$\be$}
\psline(7,0)(11,-3.5)\psline(7,-1)(11,2.5)
\rput(11.5,-2.7){$\be_1$}\rput(11.5,1.8){$\be_2$}
\psline(20,-3.5)(20,2.5)\psline(19,-2.5)(24,-3.5)\psline(19,1.5)(24,2.5)
\rput(23.5,-2.7){$\be_1$}\rput(20.7,-.5){$\be_2$}\rput(23.5,1.8){$\be_3$}
\end{pspicture}
\caption{The three possible configurations of rational curves in an ideal 
Calabi-Yau $5$-fold. The label next to each component indicates the degree.}
\label{dim5conf_fig}
\end{figure}
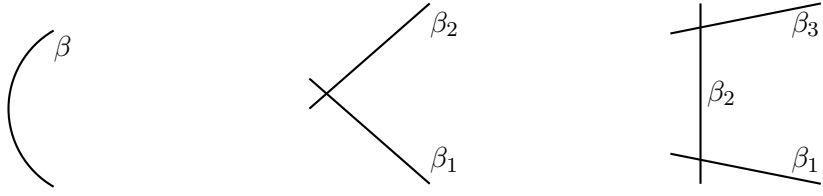

\begin{cnd}\label{g0reg_cnd}
If $u:\bP^1\!\lra\!X$ is a simple holomorphic map, $H^1\big(\bP^1,u^*TX)=0$.
\end{cnd}

\noindent
By Condition \ref{g0reg_cnd}, $\M_{0,J}^*(X,\be)$ is a nonsingular variety
of the expected dimension $2\!+\!|J|$.

\begin{cnd}\label{g0ev_cnd}
For all $\be_1,\ldots,\be_k\!\in\!H_+(X)$, finite sets $J_1,\ldots,J_k$, 
and a partition of $J_1\!\sqcup\!\ldots\!\sqcup\!J_k$ into nonempty disjoint
subsets $I_1,\ldots,I_m$, the restriction of the total evaluation 
map{\footnote{By convention, $[k]=\{1,2, \ldots, k\}$.}}
$$\ev\!:\prod_{p=1}^{k}\M_{0,J_p}^*(X,\be_p)\lra
\prod_{p=1}^{k}X^{J_p}, \qquad
\ev\big((b_p)_{p\in[k]}\big)_{(p,j)}=\ev_j(b_p)~~
\forall\,p\!\in\![k],\,j\!\in\!J_p,$$
to the open subspace of simple tuples is transverse to the diagonal
$$\big\{(x_{(p,j)})_{p\in[k],j\in J_p}\!: x_{(p,j)}\!=\!x_{(p',j')}
~\hbox{if}~(p,j),(p',j')\!\in\!I_q~\hbox{for some}~q\big\}.$$
\end{cnd}


By Condition \ref{g0ev_cnd}, $\M_{0,(\eset,\eset)}^*\big(X,(\be_1,\be_2)\big)$
and $\M_{0,\eset}^*\big(X,(\be_1,\be_2,\be_3)\big)$
are nonsingular of dimensions $1$ and $0$, respectively.
Furthermore, all simple genus $0$ maps with reducible 
domains deform to curves with nonsingular domains.
Furthermore, for all $\be\!\in\!H_+(X)$, 
the open subspace of $\ov\M_{0,J}(X,\be)$ consisting of
simple maps is nonsingular.

\begin{cnd}\label{g0der_cnd}
For all $\be\!\in\!H_+(X)$, the restriction of the bundle section $\cD_1$
to $\M_{0,1}^*(X,\be)$ is transverse to the zero set.
For all $\be_1,\be_2\!\in\!H_+(X)$, the bundle section
$$\pi_1^*\cD_0+\pi_2^*\cD_0
\in\Ga\big(\bP(\pi_1^*L_0\!\oplus\!\pi_2^*L_0)\big|_{
\M_{0,(\eset,\eset)}^*(X,(\be_1,\be_2))}, \Hom(\ga,\ev_0^*TX)\big),$$
where $\ga\!\lra\!\bP(\pi_1^*L_0\!\oplus\!\pi_2^*L_0)$ is the tautological line bundle,
is transverse to the zero set.
\end{cnd}

\noindent
By Condition \ref{g0reg_cnd} and the first part of Condition~\ref{g0der_cnd}, 
every simple holomorphic map $u\!:\bP^1\!\lra\!X$ is an immersion.
By Condition \ref{g0ev_cnd}, $u$ is injective.
Thus, every irreducible genus 0 curve $C\subset X$ is nonsingular.
The normal bundle to such a curve must split as 
$$\N=\cO(a_1)\oplus\cO(a_2)\oplus\cO(a_3)\oplus\cO(a_4)
\lra\bP^1,
\quad\hbox{with}~~~ a_i\!\in\!\Z,~\sum_{i=1}^{i=4}a_i=-2,~a_i\!\ge\!-1,$$
the last restriction follows from Condition~\ref{g0reg_cnd}.
By the first part of Condition~\ref{g0str_cnd} below, 
$a_i\!\in\!\{0,-1\}$ for all~$i$.
The second part of Condition~\ref{g0der_cnd} implies that 
every node of a reducible genus~$0$ curve in $X$ is simple.

\begin{cnd}\label{g0str_cnd}
For all $\be\!\in\!H_+(X)$, the bundle section
$$d\ev_1\in\Ga\big(\bP(T\M_{0,1}^*(X,\be)),\Hom(\ga,\ev_1^*TX)\big),$$
where $\ga\!\lra\!\bP(T\M_{0,1}^*(X,\be))$ is the tautological line bundle,
is transverse to the zero set.
For all $\be_1,\be_2\!\in\!H_+(X)$, the bundle section
$$\pi_1^*d\ev_0+\pi_2^*\cD_0\in\Ga\big(
\bP(\pi_1^*T\M_{0,\{0\}}^*(X,\be_1)\!\oplus\!\pi_2^*L_0)\big|_{
\M_{0,(\eset,\eset)}^*(X,(\be_1,\be_2))}, \Hom(\ga,\ev_0^*TX)\big),$$
where $\ga\!\lra\!\bP(\pi_1^*\!T\M_{0,\{0\}}^*(X,\be_1)\!\oplus\!\pi_2^*L_0)$ 
is the tautological line bundle, is transverse to the zero set.
\end{cnd}

\noindent
By Condition~\ref{g0str_cnd}, 
neither of the two bundle sections vanishes anywhere.
In the case of the first bundle section, the dimension of the base space and
the rank of the vector bundle both equal $5$.
On the other hand, the vanishing of the bundle section here
implies the differential of the evaluation map
$$\ev_1\!:\M_{0,1}^*(X,\be)\lra X$$
is not injective at some simple,  degree $\be$, $1$-marked map 
$[\bP^1,x_1,u]$.
Hence, the normal bundle must split as
$$\N \approx \cO(1)\oplus\cO(-1)\oplus\cO(-1)\oplus\cO(-1).$$
Therefore $d\ev_1$ is not injective at $[\bP^1,x,u]$ for all $x\!\in\!\bP^1$.
The zero set of the first bundle section in Condition~\ref{g0str_cnd}
must be at least of dimension one. So by transversality, no vanishing
is possible.


The non-vanishing of the second bundle section is clear
from transversality since the base space is of dimension 4 and
bundle is of rank 5.

\begin{lmm}
\label{g0str_lmm}
Let $X$ be an ideal Calabi-Yau $5$-fold.
If $\be\!\in\!H_+(X)$ and $J$ is a finite set, 
the space $\ov\M_{0,J}^*(X,\be)$ is nonsingular of dimension $2\!+\!|J|$
and consists of simple maps.
Furthermore, the evaluation map
$$\ev_1\!:\ov\M_{0,1}^*(X,\be)\lra X$$
is an immersion.
If $\be_1,\be_2\!\in\!H_+(X)$ and $J_1,J_2$ are finite sets, 
$\ov\M_{0,(J_1,J_2)}^*\big(X,(\be_1,\be_2)\big)$ is smooth of dimension 
$1\!+\!|J_1|\!+\!|J_2|$ and consists of simple maps.
\end{lmm}

\noindent {\em Proof.}
By Condition~\ref{g0str_cnd}, the restriction of $\ev_1$ to
the open subset 
$$\M_{0,J}^*(X,\be) \subset \ov\M_{0,J}^*(X,\be)$$
is an immersion for every $\be\!\in\!H_+(X)$.
Therefore, by the argument given in Section \ref{g1ghost_subs2},
if $$u:\Sigma \rightarrow X$$ is not simple,
then no deformation of $u$ is simple.
Hence, $u$ cannot lie in the closure of
$\M_{0,J}^*(X,\be)$.
We conclude  $\ov\M_{0,J}^*(X,\be)$ consists of simple maps
and therefore nonsingular of expected dimension.
The proof of the claim for
$\ov\M_{0,(J_1,J_2)}^*\big(X,(\be_1,\be_2)\big)$
is the same. 
\QED
\vspace{10pt}

Conditions \ref{g0reg_cnd}-\ref{g0str_cnd} can be extended to define 
an ideal Calabi-Yau $n$-fold for any~$n$.
However, Lemma~\ref{g0str_lmm}, which depends on 
the dimension counting argument in the preceding paragraph,
does not apply in dimensions $6$ and higher.
For example, if $X_8\subset \bP^7$ is the degree $8$ 
Calabi-Yau hypersurface,
$$\ov\M_{0,1}^*(X_8,1)=\ov\M_{0,1}(X_8,1)$$
certainly consists of simple maps. 
However, a computation on $G(2,8)$ shows 
the evaluation map~$\ev_1$ is not an immersion along 
$133430226944$ fibers of the forgetful morphism
$$\ov\M_{0,1}^*(X_8,1)\lra\ov\M_{0,0}^*(X_8,1).$$
A separate computation in a projective bundle over $G(3,8)$ shows 
the space of conics in $X_8$
contains $133430226944$ double lines.
In both cases the degenerate loci correspond to 
the $133430226944$ lines in $X_8$ whose normal bundle splits
as $\cO(1)\!\oplus\!\cO\!\oplus\!3\cO(-1)$, 
instead of the expected $3\cO\!\oplus\!2\cO(-1)$.
While the Calabi-Yau $6$-fold $X_8$ is not ideal, low-degree curves
in projective hypersurfaces do behave as expected.
 The appearance multiple covers as limits of simple maps 
is to be expected in dimensions~$6$ and higher,
making a full enumerative treatment more complicated
(and likely drastically so).

\subsection{Genus $0$ counts}
\label{dim5CY_subs}

We define here integer forms of the genus $0$ Gromov-Witten
invariants of 
Calabi-Yau $5$-folds by considering all possible distributions
of constraints and  $\psi$-classes between the marked points.
The $13$ relevant types of invariants are indicated in Figure
\ref{n5nums_fig}.
We  state relations motivated by ideal geometry which reduce
all 13  to genus $0$ Gromov-Witten invariants.
These relations are taken to be the definition of 13 invariants
for  arbitrary Calabi-Yau $5$-folds.

If $J$ is a finite set, $J'\!\subset\!J$, and $\be\!\in\!H_+(X)$,
let
$$f_{J,J'}\!: \ov\M_{0,J}(X,\be)\lra\ov\M_{0,J-J'}(X,\be)$$
be the forgetful map dropping the marked points indexed by the set $J'$.
If $j\!\in\!J$, let
$$\ti\psi_j=f_{J,J-j}^*\psi_j\in H^2\big(\ov\M_{0,J}(X,\be)\big),$$
where $\psi_j$ is the first chern class of the universal cotangent line bundle 
for the marked point on $\ov\M_{0,\{j\}}(X,\be)$.\\

If $X$ is an ideal Calabi-Yau $5$-fold and $\be\!\in\!H_+(X)$, 
the dimension of $\ov\M_{0,0}^*(X,\be)$ is~$2$.
There are $7$  invariants of the form
$$n_{\be}(\ti\psi^a\mu_1,\mu_2,\ldots,\mu_k)
=\int_{\ov\M_{0,k}^*(X,\be)}\ti\psi_1^a\prod_{j=1}^{k}\ev_j^*\mu_j,
\qquad a\!\ge\!0,~\mu_j\!\in\!H^{2*}(X),$$
which we require:
\begin{enumerate}[label=(1\Alph*)]
\item\label{d5n1p1psi0}
$n_{\be}(\mu)$ where $\mu\!\in\!H^6(X)$ counting curves through $\mu$,
\item\label{d5n1p2psi0}
$n_{\be}(\mu_1,\mu_2)$ where $\mu_1,\mu_2\!\in\!H^4(X)$ counting
curves through $\mu_1$ and $\mu_2$,
\item\label{d5n1p1psi1} $n_{\be}(\ti\psi\mu)$  where $\mu\!\in\!H^4(X)$,
\item\label{d5n1p2psi1} 
$n_{\be}(\ti\psi\mu_1,\mu_2)$ 
where $\mu_1\!\in\!H^2(X)$ and $\mu_2\!\in\!H^4(X)$,
\item\label{d5n1p1psi2} $n_{\be}(\ti\psi^2\mu)$ where $\mu\!\in\!H^2(X)$,
\item\label{d5n1p2psi2} 
$n_{\be}(\ti\psi^2,\mu)$ where $\mu\!\in\!H^4(X)$,
\item\label{d5n1p1psi3} $n_{\be}(\ti\psi^3)$.
\end{enumerate}
Let $\ov\cM_{\be}$ denote the unpointed space $\ov\M_{0,0}^*(X,\be)$.
We will  need the Chern number
\begin{enumerate}[label=(1H)]
\item\label{d5n1euler} 
$\gamma_1(\be)=\int_{\ov\cM_\be} 
\big(c_1^2(\ov\cM_{\be})-c_2(\ov\cM_{\be})\big).$
\end{enumerate}

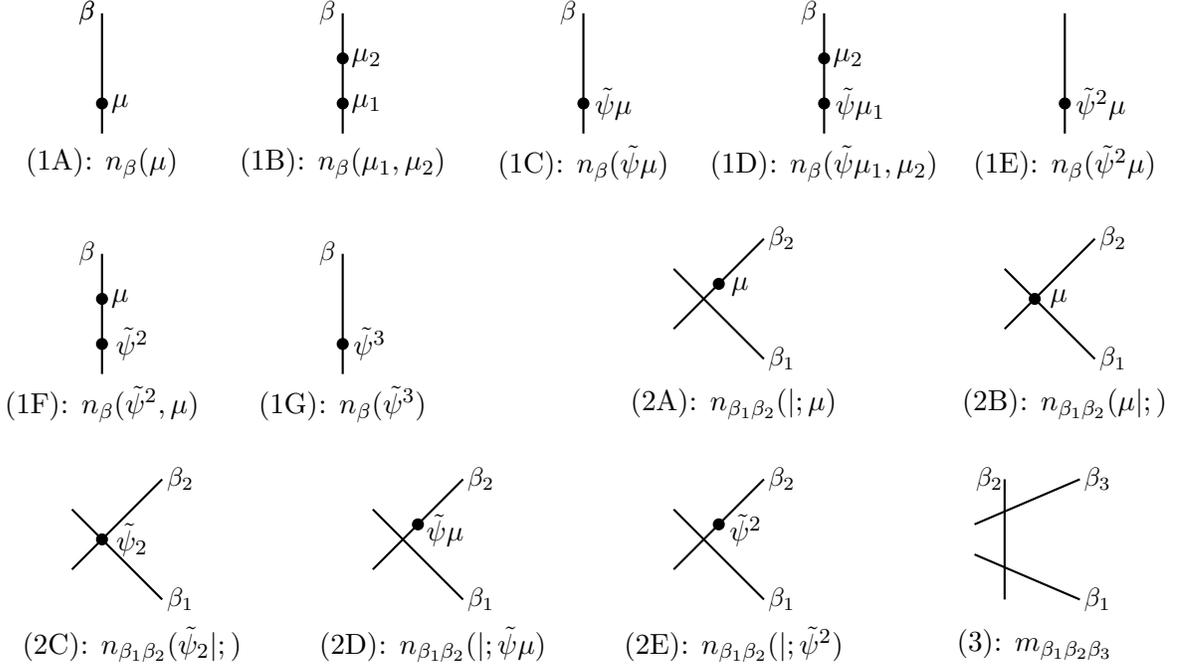
\begin{figure}
\begin{pspicture}(-3,-8)(10,1.3)
\psset{unit=.4cm}
\psline(-2,-2)(-2,2)\rput(-2.5,2){\sm{$\be$}}
\pscircle*(-2,-1){.2}\rput(-1.4,-1){$\mu$}
\rput(-2,-3){\ref{d5n1p1psi0}: $n_{\be}(\mu)$}
\psline(6,-2)(6,2)\rput(5.5,2){\sm{$\be$}}
\pscircle*(6,-1){.2}\rput(6.8,-1){$\mu_1$}
\pscircle*(6,.5){.2}\rput(6.8,.5){$\mu_2$}
\rput(6,-3){\ref{d5n1p2psi0}: $n_{\be}(\mu_1,\mu_2)$}
\psline(14,-2)(14,2)\rput(13.5,2){\sm{$\be$}}
\pscircle*(14,-1){.2}\rput(15,-1){$\ti\psi\mu$}
\rput(14,-3){\ref{d5n1p1psi1}: $n_{\be}(\ti\psi\mu)$}
\psline(22,-2)(22,2)\rput(21.5,2){\sm{$\be$}}
\pscircle*(22,-1){.2}\rput(23.2,-1){$\ti\psi\mu_1$}
\pscircle*(22,.5){.2}\rput(22.8,.5){$\mu_2$}
\rput(22,-3){\ref{d5n1p2psi1}: $n_{\be}(\ti\psi\mu_1,\mu_2)$}
\psline(30,-2)(30,2)\rput(-2.5,2){\sm{$\be$}}
\pscircle*(30,-1){.2}\rput(31.2,-1){$\ti\psi^2\mu$}
\rput(30,-3){\ref{d5n1p1psi2}: $n_{\be}(\ti\psi^2\mu)$}
\psline(-2,-10)(-2,-6)\rput(-2.5,-6){\sm{$\be$}}
\pscircle*(-2,-9){.2}\rput(-1,-9){$\ti\psi^2$}
\pscircle*(-2,-7.5){.2}\rput(-1.4,-7.5){$\mu$}
\rput(-2,-11){\ref{d5n1p2psi2}: $n_{\be}(\ti\psi^2,\mu)$}
\psline(6,-10)(6,-6)\rput(5.5,-6){\sm{$\be$}}
\pscircle*(6,-9){.2}\rput(6.9,-9){$\ti\psi^3$}
\rput(6,-11){\ref{d5n1p1psi3}: $n_{\be}(\ti\psi^3)$}
\psline(17,-8.5)(20,-5.5)\rput(20.6,-5.5){\sm{$\be_2$}}
\psline(17,-6.5)(20,-9.5)\rput(20.6,-9.5){\sm{$\be_1$}}
\pscircle*(18.5,-7){.2}\rput(19.2,-7.1){$\mu$}
\rput(19,-11){\ref{d5n2p1psi0}: $n_{\be_1\be_2}(|;\mu)$}
\psline(28,-8.5)(31,-5.5)\rput(31.6,-5.5){\sm{$\be_2$}}
\psline(28,-6.5)(31,-9.5)\rput(31.6,-9.5){\sm{$\be_1$}}
\pscircle*(29,-7.5){.2}\rput(29.8,-7.5){$\mu$}
\rput(30,-11){\ref{d5n2p0psi0}: $n_{\be_1\be_2}(\mu|;)$}
\psline(-3,-16.5)(0,-13.5)\rput(0.6,-13.5){\sm{$\be_2$}}
\psline(-3,-14.5)(0,-17.5)\rput(0.6,-17.5){\sm{$\be_1$}}
\pscircle*(-2,-15.5){.2}\rput(-1,-15.5){$\ti\psi_2$}
\rput(-1,-19){\ref{d5n2p0psi1}: $n_{\be_1\be_2}(\ti\psi_2|;)$}
\psline(7,-16.5)(10,-13.5)\rput(10.6,-13.5){\sm{$\be_2$}}
\psline(7,-14.5)(10,-17.5)\rput(10.6,-17.5){\sm{$\be_1$}}
\pscircle*(8.5,-15){.2}\rput(9.4,-15.2){$\ti\psi\mu$}
\rput(9,-19){\ref{d5n2p1psi1}: $n_{\be_1\be_2}(|;\ti\psi\mu)$}
\psline(17,-16.5)(20,-13.5)\rput(20.6,-13.5){\sm{$\be_2$}}
\psline(17,-14.5)(20,-17.5)\rput(20.6,-17.5){\sm{$\be_1$}}
\pscircle*(18.5,-15){.2}\rput(19.4,-15.2){$\ti\psi^2$}
\rput(19,-19){\ref{d5n2p1psi2}: $n_{\be_1\be_2}(|;\ti\psi^2)$}
\psline(28,-17.5)(28,-13.5)\rput(27.5,-13.5){\sm{$\be_2$}}
\psline(27,-15)(30.5,-13.5)\rput(31.1,-13.5){\sm{$\be_3$}}
\psline(27,-16)(30.5,-17.5)\rput(31.1,-17.5){\sm{$\be_1$}}
\rput(29,-19){\ref{d5n3}: $m_{\be_1\be_2\be_3}$}
\end{pspicture}
\caption{Counts for Calabi-Yau $5$-folds}
\label{n5nums_fig}
\end{figure}

There are $5$ types of relevant counts of connected 2-component curves
which we require,
\begin{equation*}\begin{split}
&n_{\be_1\be_2}\big(\ti\psi_1^{a_1}\ti\psi_2^{a_2}\mu_0|
\ti\psi^{b_1}\mu_{1,1},\mu_{1,2},\ldots,\mu_{1,k_1};
\ti\psi^{b_2}\mu_{2,1},\mu_{2,2},\ldots,\mu_{2,k_2}\big)\\
&\qquad\qquad
=\int_{\ov\M_{0,([k_1],[k_2])}^*(X,(\be_1,\be_2))}
\pi_1^*\bigg(\ti\psi_0^{a_1}\ti\psi_1^{b_1}\ev_0^*\mu_0
\prod_{j=1}^{k_1}\!\ev_j^*\mu_{1,j}\bigg) 
\pi_2^*\bigg(\ti\psi_0^{a_2}\ti\psi_1^{b_2}\prod_{j=1}^{k_2}\!\ev_j^*\mu_{2,j}\bigg),
\end{split}\end{equation*}
where $\pi_1,\pi_2$ are the component projection maps as in~\e_ref{pi12_e},
$a_i,b_i\!\ge\!0$, and $\mu_0,\mu_{i,j}\!\in\!H^{2*}(X)$. 
The 
$5$ types are represented by the following counts of $(\be_1,\be_2)$-curves: 
\begin{enumerate}[label=(2\Alph*)]
\item\label{d5n2p1psi0}
$n_{\be_1\be_2}(|;\mu)$ where $\mu\!\in\!H^4(X)$,
\item\label{d5n2p0psi0}
$n_{\be_1\be_2}(\mu|;)$ where $\mu\!\in\!H^2(X)$,
\item\label{d5n2p0psi1}
$n_{\be_1\be_2}(\ti\psi_2|;)$,
\item\label{d5n2p1psi1}
$n_{\be_1\be_2}(|;\ti\psi\mu)$ where $\mu\!\in\!H^2(X)$,
\item\label{d5n2p1psi2}
$n_{\be_1\be_2}(|;\ti\psi^2)$.
\end{enumerate}

Finally, we denote the cardinality of the compact $0$-dimensional space
$\M_{0,\eset}^*(\be_1,\be_2,\be_3)$ for triples 
$\be_1,\be_2,\be_3\!\in\!H_+(X)$
by $m_{\be_1\be_2\be_3}$:
\begin{enumerate}[label=(3)]
\item\label{d5n3}
$m_{\be_1\be_2\be_3}$ is the number of connected 3-component curves
of tridegree $\be_1,\be_2,\be_3$.
\end{enumerate}

\begin{figure}
\begin{pspicture}(-3.5,-1.4)(10,1.6)
\psset{unit=.4cm}
\rput(-5,0){\Lg{$\psi_1=$}}
\psline(-2,-3)(-2,2)\rput(-2.5,2){\sm{$\be$}}
\pscircle*(-2,-2){.2}\rput(-1.5,-2){\sm{$1$}}
\psline(-3,0)(0,3)\rput(.5,3.1){\sm{$0$}}
\pscircle*(-1.33,1.66){.2}\rput(-1.33,1){\sm{$2$}}
\pscircle*(-.66,2.33){.2}\rput(-.66,1.7){\sm{$3$}}
\rput(3,-1){\Lg{$+\sum\limits_{\underset{\be_1,\be_2\in H_+(X)}{\be_1+\be_2=\be}}$}}
\psline(7,-3)(7,2)\rput(6.5,2){\sm{$\be_2$}}
\pscircle*(7,-2){.2}\rput(7.5,-2){\sm{$1$}}
\psline(6,0)(9,3)\rput(9.5,3.1){\sm{$\be_1$}}
\pscircle*(7.67,1.66){.2}\rput(7.67,1){\sm{$2$}}
\pscircle*(8.34,2.33){.2}\rput(8.34,1.7){\sm{$3$}}
\rput(13.7,0){\Lg{$\ti\psi_1=\psi_1-$}}
\psline(18,-3)(18,2)\rput(17.5,2){\sm{$0$}}
\pscircle*(18,-2){.2}\rput(18.5,-2){\sm{$1$}}
\psline(17,0)(20,3)\rput(20.5,3.1){\sm{$\be$}}
\pscircle*(18,-.5){.2}\rput(18.5,-.5){\sm{$2$}}
\pscircle*(19,2){.2}\rput(19,1.36){\sm{$3$}}
\rput(20.4,-.2){\Lg{$-$}}
\psline(23,-3)(23,2)\rput(22.5,2){\sm{$0$}}
\pscircle*(23,-2){.2}\rput(23.5,-2){\sm{$1$}}
\psline(22,0)(25,3)\rput(25.5,3.1){\sm{$\be$}}
\pscircle*(23,-.5){.2}\rput(23.5,-.5){\sm{$3$}}
\pscircle*(24,2){.2}\rput(24,1.36){\sm{$2$}}
\rput(25.4,-.2){\Lg{$-$}}
\psline(28,-3)(28,2)\rput(27.5,2){\sm{$0$}}
\pscircle*(28,-2.25){.2}\rput(28.5,-2.25){\sm{$1$}}
\psline(27,0)(30,3)\rput(30.5,3.1){\sm{$\be$}}
\pscircle*(28,-.25){.2}\rput(28.5,-.25){\sm{$3$}}
\pscircle*(28,-1.25){.2}\rput(28.5,-1.25){\sm{$2$}}
\end{pspicture}
\caption{Relations for $\psi_1$ and $\ti\psi_1$ on $\ov\M_{0,3}^*(X,\be)$.
Each curve represents the divisor in $\ov\M_{0,3}^*(X,\be)$ whose
general element has the domain and the degree distribution specified by the curve.}
\label{TRC_fig}
\end{figure}
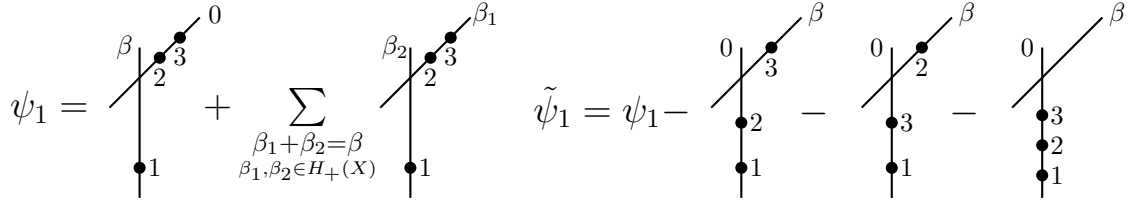

The numbers \ref{d5n1p1psi0} and \ref{d5n1p2psi0} are determined from
$1$- and $2$-pointed Gromov-Witten invariants via \e_ref{g0BPS_e}.
The tautological recursion relation for $\psi_1$ can be used 
to express $\ti\psi_1$ in terms of boundary divisors on 
$\ov\M_{0,3}^*(X,\be)$, see Figure \ref{TRC_fig}.
The divisor relation then gives rise to the relations between the invariants 
\ref{d5n1p1psi1}-\ref{d5n1p1psi3} indicated in Figure \ref{tipsi_fig},
see also Section 3 in \cite{g0inter}.
We now describe these relations formally.
If $H$ is a divisor on~$X$ and $H_{\be}\!=\!(H,\be)$, then
\BE{n1psired_e}\begin{split}
H_{\be}^2\,n_{\be}(\ti\psi\mu)&=
n_{\be}(\mu,H^2)-2H_{\be}\,n_{\be}(H\mu)+
\sum_{\be_1+\be_2=\be}\!\!\!\!\!H_{\be_1}^2 n_{\be_1\be_2}(|;\mu),\\
H_{\be}^2\,n_{\be}(\ti\psi\mu_1,\mu_2)&=
(\mu_1,\be)\,n_{\be}(\mu_2,H^2)-2H_{\be}\,n_{\be}(H\mu_1,\mu_2)\\
&\qquad\qquad+
\sum_{\be_1+\be_2=\be}\!\!\!\!\!
\big((\mu_1,\be_1)H_{\be_2}^2\!+\!(\mu_1,\be_2)H_{\be_1}^2\big)
n_{\be_1\be_2}(|;\mu_2),\\
H_{\be}^2\,n_{\be}(\ti\psi^2\mu)&=
n_{\be}(\ti\psi\mu,H^2)-2H_{\be}\,n_{\be}(\ti\psi H\mu)+
\sum_{\be_1+\be_2=\be}\!\!\!\!\!H_{\be_1}^2
\big(n_{\be_1\be_2}(|;\ti\psi\mu)\!+\!n_{\be_1\be_2}(\mu|;)\big),\\
n_{\be}(\ti\psi^2,\mu)&=-\sum_{\be_1+\be_2=\be}\!\!\!\!
\!n_{\be_1\be_2}(|;\mu),\\
H_{\be}^2\,n_{\be}(\ti\psi^3)&
=n_{\be}(\ti\psi^2,H^2)-2H_{\be}\,n_{\be}(\ti\psi^2H)
+\sum_{\be_1+\be_2=\be}\!\!\!\!\!H_{\be_1}^2\big(
n_{\be_1\be_2}(|;\ti\psi^2)+n_{\be_1\be_2}(\ti\psi_2|;)\big),
\end{split}\EE
the fourth identity above is obtained by applying the relation of 
Figure \ref{tipsi_fig} twice.
We can similarly remove $\psi$-classes from 2-component curves:
\BE{n2psired_e}\begin{split}
H_{\be_2}^2\,n_{\be_1\be_2}(\ti\psi_2|;)&=
n_{\be_1\be_2}(|;H^2)-2H_{\be_2}\,n_{\be_1\be_2}(H|;)
+\sum_{\be+\be'=\be_2}\!\!\!H_{\be}^2\,m_{\be_1\be'\be},\\
H_{\be_2}^2\,n_{\be_1\be_2}(|;\ti\psi\mu)&=
(\mu,H)\,n_{\be_1\be_2}(|;H^2)-2H_{\be_2}\,n_{\be_1\be_2}(|;H\mu)\\
&\qquad\qquad
+\sum_{\be+\be'=\be_2}\!\!\!
\big((\mu,\be)H_{\be'}^2\!+\!(\mu,\be')H_{\be}^2\big)m_{\be_1\be'\be},\\
n_{\be_1\be_2}(|;\ti\psi^2)&=-\sum_{\be+\be'=\be_2}\!\!\!m_{\be_1\be'\be}\,,
\end{split}\EE
the last identity above is obtained by applying the relation of 
Figure \ref{tipsi_fig} twice.
On the other hand, by \e_ref{cF_e} and some manipulation,
\BE{etanum_e}\begin{split}
\gamma_1(\be)&=\frac{1}{2}\Big(n_{\be}\big(c_3(X)\big)
+n_{\be}\big(\ti\psi c_2(X)\big)+n_{\be}\big(\ti\psi^3\big)
+n_{\be}\big(c_2(X),c_2(X)\big)+4\,n_{\be}\big(\ti\psi^2,c_2(X)\big)\Big)\\
&\qquad\qquad
-\sum_{\be_1+\be_2=\be}\Big(2\,n_{\be_1\be_2}(|;\ti\psi^2)
+\frac{5}{2}\,n_{\be_1\be_2}(\ti\psi_2|;)\Big).
\end{split}\EE

\begin{figure}
\begin{pspicture}(-2.7,-2.2)(10,1.6)
\psset{unit=.4cm}
\rput(-3,0){$H_{\be}^2$}
\psline(-2,-3)(-2,2)\rput(-2,2.6){\sm{$\be$}}
\pscircle*(-2,-2){.2}\rput(-.7,-1.8){\sm{$\ti\psi^c\mu_e$}}
\rput(-2.7,-1.8){\sm{$e$}}
\rput(.6,0){\Lg{$=$}}
\psline(3,-3)(3,2)\rput(3,2.6){\sm{$\be$}}
\pscircle*(3,-2){.2}\rput(4.8,-1.8){\sm{$\ti\psi^{c-1}\mu_e$}}
\rput(2.3,-1.8){\sm{$e$}}
\pscircle*(3,.5){.2}\rput(4,.6){\sm{$H^2$}}
\rput(2.3,.5){\sm{$0$}}
\rput(6.5,0){\Lg{$-$}}\rput(8.2,-.1){$2H_{\be}$}
\psline(9.5,-3)(9.5,2)\rput(9.5,2.6){\sm{$\be$}}
\pscircle*(9.5,-2){.2}\rput(11.8,-1.8){\sm{$\ti\psi^{c-1}H\mu_e$}}
\rput(8.8,-1.8){\sm{$e$}}
\rput(16.5,-1){\Lg{$+
\sum\limits_{\underset{\be_1,\be_2\in H_+(X)}{\be_1+\be_2=\be}}$}}
\rput(20.3,0){$H_{\be_1}^2$}\rput(22,0){\Hg{$\Bigg\{$}}
\psline(24,-3)(24,2)\rput(24,2.6){\sm{$\be_2$}}
\pscircle*(24,-2){.2}\rput(25.8,-1.8){\sm{$\ti\psi^{c-1}\mu_e$}}
\rput(23.5,-1.8){\sm{$e$}}
\psline(23,0)(25.5,2.5)\rput(26,2.6){\sm{$\be_1$}}
\rput(26.5,0){\Lg{$+$}}
\psline(29,-3)(29,2)\rput(29,2.6){\sm{$\be_2$}}
\psline(28,0)(30.5,2.5)\rput(31,2.6){\sm{$\be_1$}}
\pscircle*(29,1){.2}\rput(30.8,.5){\sm{$\ti\psi_2^{c-2}\mu_{e}$}}
\rput(28.3,1.3){\sm{$e$}}\rput(32.5,0){\Hg{$\Bigg\}$}}\rput(15,-5){\sm{$H\subset X$ divisor,
$H_{\be}=H\cdot\be,~H_{\be_1}=H\cdot\be_1$}}
\end{pspicture}
\caption{Reducing the power of $\ti\psi$ at marked point $e$ 
in the absence of $\psi$-classes at other marked points.}
\label{tipsi_fig}
\end{figure}

The meeting numbers \ref{d5n2p1psi0}, \ref{d5n2p0psi0}, and~\ref{d5n3}
are computed via degree reducing recursions analogous to 
Rules (i)-(iv) of Section 0.3 of \cite{KlP} for the $4$-dimensional case.
Let 
$$\{\om_1,\ldots,\om_N\},\{\om_1^{\#}\,\ldots,\om_N^{\#}\}
\subset H^4(X)\oplus H^6(X) $$
be dual bases normalized so that 
$$\PD_{X^2}\De_X-\sum_{l=1}^{N}\om_l\!\times\!\om_l^{\#}
\in \bigoplus_{k=0,1,4,5}\!\!\!\!\!
H^{2k}(X)\!\otimes\!H^{2(5-k)}(X) 
\oplus H^{odd}(X)\!\otimes\!H^{odd}(X),$$
where $\De_X\!\subset\!X^2$ is the diagonal.
Then,
\BE{m2red_e1}\begin{split}
n_{\be_1\be_2}(|;\mu)&=
\sum_{l=1}^{N}\! n_{\be_1}(\om_l)\,n_{\be_2}(\om_l^{\#},\mu)
+\begin{cases}
n_{\be_1,\be_2-\be_1}(|;\mu)+n_{\be_2-\be_1,\be_1}(|;\mu),&
\hbox{if}~\be_2\!>\!\be_1,\\
n_{\be_1-\be_2,\be_2}(|;\mu),& \hbox{if}~\be_2\!<\!\be_1,\\
n_{\be_1}(c_2(X),\mu)+
2n_{\be_1}(\ti\psi^2,\mu),& \hbox{if}~\be_2\!=\!\be_1.
\end{cases}\end{split}\EE
In light of the fourth identity in \e_ref{n1psired_e},
the relation differs from the $4$-dimensional case only 
by the expected adjustment for the constraint $\mu$.

The corresponding recursions for the numbers \ref{d5n2p0psi0} and~\ref{d5n3}
are more complicated.
For classes $\be_1,\be_2\!\in\!H_+(X)$, let
\BE{eta2dfn_e}
\gamma_2(\be_1,\be_2)=n_{\be_1\be_2}\big(|;c_2(X)\big)
+2\,n_{\be_1\be_2}(|;\ti\psi^2)+
n_{\be_1\be_2}(\ti\psi_2|;)+n_{\be_2\be_1}(\ti\psi_2|;).
\EE
For $\mu\!\in\!H^2(X)$, we define
\BE{n2corr_e0}
\tC_{\be_1\be_2}(\mu)=\begin{cases}
\begin{aligned}
&n_{\be_2-\be_1,\be_1}(|;\ti\psi\mu)+n_{\be_2-\be_1,\be_1}(\mu|;)\\
&\qquad +(\mu,\be_1)\big(\gamma_2(\be_2\!-\!\be_1,\be_1)+\frac{1}{2}
\sum\limits_{\be+\be'=\be_2-\be_1}\!\!\!\!\!m_{\be\be_1\be'}\big),
\end{aligned}
&\hbox{if}~\be_2\!>\!\be_1;\\
\tC_{\be_2\be_1}(\mu), &\hbox{if}~\be_2\!<\!\be_1;\\
\begin{aligned}\\
&n_{\be_1}(c_2(X)\mu)+n_{\be_1}(\ti\psi^2\mu)+n_{\be_1}(c_2(X),\ti\psi\mu)
+(\mu,\be_1)\gamma_1(\be_1)\\
&\quad\qquad-\sum_{\be+\be'=\be_2}\!\!(2n_{\be\be'}(|;\ti\psi\mu)+
\frac{5}{2}n_{\be\be'}(\mu|;)\big),
\end{aligned}&\hbox{if}~\be_1\!=\!\be_2.
\end{cases}
\EE
For $\be_1,\be_2,\be_3\!\in\!H_+(X)$, let 
\BE{n3corr_e}\begin{split}
\tC^{(1)}_{\be_1\be_2\be_3}&=
\begin{cases}
m_{\be_3-\be_1,\be_1,\be_2},&\hbox{if}~\be_3\!>\!\be_1;\\
m_{\be_1-\be_3,\be_3,\be_2},&\hbox{if}~\be_3\!<\!\be_1;\\
\gamma_2(\be_2,\be_1),&\hbox{if}~\be_3\!=\!\be_1;
\end{cases}\\
\tC^{(2)}_{\be_1\be_2\be_3}&=-
\begin{cases}
m_{\be_1,\be_2,\be_3-\be_2},&\hbox{if}~\be_3\!>\!\be_2;\\
m_{\be_1,\be_3,\be_2-\be_3}\!+\!m_{\be_1,\be_2-\be_3,\be_3},
&\hbox{if}~\be_3\!<\!\be_2;\\
n_{\be_1\be_2}\big(|;c_2(X)\big)+2\,n_{\be_1\be_2}(|;\ti\psi^2),
&\hbox{if}~\be_3\!=\!\be_2;
\end{cases}\\
\tC^{(12)}_{\be_1\be_2\be_3}&=
-\begin{cases}
m_{\be_3-\be_1-\be_2,\be_1,\be_2},&\hbox{if}~\be_3\!>\!\be_1\!+\!\be_2;\\
m_{\be_1\!+\!\be_2\!-\!\be_3,\be_3-\be_2,\be_2},
&\hbox{if}~\be_2\!<\!\be_3\!<\!\be_1\!+\!\be_2;\\
\gamma_2(\be_2,\be_1), &\hbox{if}~\be_3\!=\!\be_1\!+\!\be_2;\\
0,&\hbox{otherwise}.
\end{cases}
\end{split}\EE
Then,
\begin{alignat}{1}
\label{rec_e2a}
n_{\be_1\be_2}(\mu|;)&=
\sum_{l=1}^{N}\! n_{\be_1}(\om_l\mu)\,n_{\be_2}(\om_l^{\#})
-\sum_{\be<\be_1,\be_2}\!\!\!(\mu,\be)\,m_{\be_1-\be,\be,\be_2-\be}
-\tC_{\be_1\be_2}(\mu),\\
\label{rec_e2b}
m_{\be_1\be_2\be_3}&=\sum_{l=1}^{N}\! n_{\be_1\be_2}(|;\om_l)\,n_{\be_3}(\om_l^{\#})
-\tC_{\be_1\be_2\be_3}^{(1)}-\tC_{\be_1\be_2\be_3}^{(2)}
-\tC_{\be_1\be_2\be_3}^{(12)} \,.
\end{alignat}
A few low degree 2-component meeting numbers for a degree~$7$
hypersurface in~$\bP^6$ are given in Table \ref{g0_table2}.
The number $n_{1,1}(H|;)$ can be confirmed via a Schubert computation
similar to Section 3 in \cite{Ka}.\\

\begin{table}
\begin{center}
\begin{tabular}{|r|r|r|}
\hline
$n_{d_1d_2}(H|;)$& $d_2=1$& $d_2=2$\\
\hline
$d_1=1$& 145366465734& 17628837973096812\\
$2$& 17628837973096812& 2134616449608028257452\\
$3$&  4403307962301366086458& 533112594803936499402982169\\
\hline
\end{tabular}
\end{center}
\caption{Meeting invariants $n_{d_1d_2}(H|;)$ for a degree $7$ hypersurface 
in $\bP^6$ counting the virtual number of
 $(d_1,d_2)$-curves with node on a fixed hyperplane.}
\label{g0_table2}
\end{table}

\noindent

Configurations of rational curves 
in a Calabi-Yau $n$-fold for can be studied for any $n$.
If $n\!\ge\!6$, such configurations include curves with 
non-simple nodes (several components sharing a node).
While describing such curves is just notationally involved,
specifying degree reducing recursions for them (following the approach
of Section \ref{recpf_subs} below) presents new difficulties.
In particular, curves with unbalanced splittings of the normal bundle 
will effect excess contributions via the loci of non-simple tuples
of maps in the closures of simple tuples of maps, 
see the end of Section \ref{g0config_subs}.
Thus, separate counts must be set up for such curves, and 
their multiple-cover contributions to 
the appropriate topological intersection numbers 
(represented by the first terms on the right-hand side
of \e_ref{m2red_e1}, \e_ref{rec_e2a}, and \e_ref{rec_e2b}) must be determined.

\subsection{Justification of degree reducing recursions}
\label{recpf_subs}

\subsubsection{Overview}
Each curve $\cC$ of type \ref{d5n2p1psi0}, \ref{d5n2p0psi0}, and \ref{d5n3}
determines a pair $(\bar\cC,\cC^*)$ of curves, where $\cC^*$ is the last
component of $\cC$ and $\bar\cC$ consists of the remaining component(s)
of $\cC$. 
The curve $\bar\cC$ has 1 component 
in the first two cases and 2 components in the last case.
The curves $\bar\cC$ and $\cC^*$ carry marking 
$x_e\!\in\!\bar\cC$ and $y_e\!\in\!\cC^*$ satisfying $x_e\!=\!y_e$.
We denote by $\ov\cM$ and $\cM^*$ the corresponding compactified spaces 
of curves/maps:
\begin{equation*}
\begin{array}{lcc}
& \ov\cM&   \cM^*\\
\hbox{Case~\ref{d5n2p1psi0}:}\qquad& \ov\M_{0,\{e\}}^*(X,\be_1)&
\big\{\phi\!\in\!\ov\M_{0,\{e\}}^*(X,\be_2)\!: (\Im\,\phi)
\!\cap\!\mu\!\neq\!\eset\big\},\\
\hbox{Case~\ref{d5n2p0psi0}:}\qquad&
\big\{\phi\!\in\!\ov\M_{0,\{e\}}^*(X,\be_1)\!: \ev_e(\phi)\!\in\!\mu\big\}&
\ov\M_{0,\{e\}}^*(X,\be_2),\\
\hbox{Case~\ref{d5n3}:}\qquad&
\ov\M_{0,(\eset,\{e\})}^*(X,(\be_1,\be_2))& \ov\M_{0,\{e\}}^*(X,\be_3),
\end{array}\end{equation*} 
where $\mu$ above denotes a generic representative for the Poincare dual
of $\mu\!\in\!H^*(X)$.
The evaluation map
$$\ev_{e,e}\!:\ov\cM\times\cM^*\lra X\!\times\!X, \quad 
\big((\bar\cC,x_e),(\cC^*,y_e)\big) \lra (x_e,y_e),$$
is then a cycle of (complex) dimension~$5$.
The relevant meeting number is the cardinality of the subset of
$$\cZ= \ev_{e,e}^{-1}(\De_X)
=\big\{ \big((\bar\cC,x_e),(\cC^*,y_e)\big)\!\in\!\ov\cM\!\times\!\cM^*\!:
x_e\!=\!y_e\big\}$$
consisting of simple pairs of maps.

The homological intersection number of the cycle $\ev_{e,e}$ with 
the class of the diagonal $\De_X\!\subset\!X^2$ in $X^2$ is 
given by the diagonal-splitting term on the right-hand side
of \e_ref{m2red_e1}, \e_ref{rec_e2a}, and \e_ref{rec_e2b}.
The homological intersection
is the number of points, counted with sign, in the preimage of $\De_X$
under a small deformation of the map $\ev_{e,e}$.
All such points must lie near $\cZ$. 
The points of $\cZ$ at which $\ev_{e,e}$ is transverse to $\De_X$
contribute $1$ each to the homology intersection.
These points include all tuples as above such that the curves $\bar\cC$ and
 $\cC^*$  do not have any components in common.
Thus, the relevant meeting number is the diagonal-splitting term 
in \e_ref{m2red_e1}, \e_ref{rec_e2a}, and \e_ref{rec_e2b} minus the contribution to 
the homology intersection number  of $\ev_{e,e}$ with $\De_X$ from 
the subset $\cZ'$  of $\cZ$ consisting of tuples as above 
such that $\bar\cC$ and $\cC^*$ have at least one component in common.
In the rest of this subsection, 
we determine these tuples and their excess 
contributions.\footnote{As 
in the $4$-dimensional case considered in~\cite{KlP},
all contributions in case (2A) are degenerate contributions arising
from loci of dimensions 1 and 2.
However, in cases \ref{d5n2p0psi0} and \ref{d5n3},
$\cZ'$ includes regular points  with respect to the evaluation condition which
are isolated and nondegenerate.}

If $X$ is an ideal Calabi-Yau $5$-fold and $\be\!\in\!H_+(X)$, the space 
$$\ov\cM_{\be,1}=\ov\M_{0,1}^*(X,\be)$$
of simple maps to $X$ of degree $\be$ with $1$ marking 
is nonsingular of dimension $3$, and the evaluation map
$$\ev\!:\ov\cM_{\be,1}\lra X$$
is an immersion, see Lemma~\ref{g0str_lmm}.
We denote by $T_{\be}$ the tangent bundle of $\ov\cM_{\be,1}$
and by $\N_{\be}$ the normal bundle to the immersion~$\ev$.
Let $\N_{\De}\!\lra\!\De$ be the normal bundle to the diagonal in~$X^2$.
If $\cC\!\subset\!X$ is a curve, let $|\cC|$ denote the number of
irreducible components of~$\cC$.

\subsubsection{Chern classes}
\label{cherh_sssub}

\noindent 
Let $X$ be an ideal Calabi-Yau $5$-fold, and let
$\be\!\in\!H_+(X)$.
We relate here the Chern classes of the normal bundle
$\N_{\be}$ to the immersion
$$\ev_1\!: \ov\cM_{\be,1}\lra X$$
to meeting numbers. 
Denote by
\BE{beforgdfn_e} f\!:\ov\cM_{\be,1}\lra\ov\cM_{\be}\EE
the forgetful map to  the nonsingular 2-dimensional
moduli space
$\ov\cM_{\be}\!=\!\ov\M_{0,0}^*(X,\be)$.

Using the bundle homomorphism $df\!:T\ov\cM_{\be,1}\!\lra\!f^*T\ov\cM_{\be}$ 
over $\ov\cM_{\be,1}$, we obtain
\begin{equation}\label{cS_e}\begin{split}
c_1(T_{\be})&=-\psi+f^*c_1\big(\ov\cM_{\be}\big),\\
c_2(T_{\be})&=\De-\psi\,f^*c_1\big(\ov\cM_{\be}\big)
+f^*c_2\big(\ov\cM_{\be}\big),
\end{split}\end{equation}
where $\psi$ is the first chern class of the cotangent line bundle on 
$\ov\cM_{\be,1}$ viewed as a $1$-pointed moduli space
and $\De\!\subset\!\ov\cM_{\be,1}$ is the locus of 
singular points of~$f$ (points at which $df$ is not surjective).
On the other hand, since $c_1(X)\!=\!0$, 
\BE{tanvsnorm_e}\begin{split}
c_1\big(\N_{\be}\big)&=-c_1\big(T_{\be}\big),\\
c_2\big(\N_{\be}\big) &= 
\ev^*c_2(X)+c_1^2\big(T_{\be}\big)-c_2\big(T_{\be}\big).
\end{split}\EE
Combining  \e_ref{cS_e} and \e_ref{tanvsnorm_e}, we find 
\BE{cN_e}\begin{split}
c_1(\N_{\be}) &= \psi-f^*c_1\big(\ov\cM_{\be}\big),\\
c_2(\N_{\be})&=\ev^*c_2(X)+\psi^2-\De-\psi\,f^*c_1\big(\ov\cM_{\be}\big)
+f^*\big(c_1^2(\ov\cM_{\be})-c_2(\ov\cM_{\be})\big).
\end{split}\EE

If $\be_1\!+\!\be_2\!=\!\be$ and $\be_1\!\neq\!\be_2$, 
let $D_{\be_1,\be_2}\!\subset\!\ov\cM_{\be}$ be the closure of 
the locus consisting of $\be$-curves split into a $\be_1$-curve and a $\be_2$-curve.
If $2\be_1\!=\!\be$, let $D_{\be_1\be_1}\!\subset\!\ov\cM_{\be}$ be twice
the closure of the locus of consisting of $\be$-curves split into two $\be_1$-curves.
In particular,
$$f_*\De=\frac{1}{2}\sum_{\underset{\be_1,\be_2\in H_+(X)}{\be_1+\be_2=\be}}
\!\!\!\!\!\!\!\!\!D_{\be_1,\be_2}\,.$$
Denote by
$(\psi_1\!+\!\psi_2)D_{\be_1,\be_2}\in H^4(\ov\cM_{\be})$
the class obtained by capping $\De$ with the first chern class of
the cotangent line bundle at the chosen node for each of the two curves.
From a Grothendieck-Riemann-Roch computation
applied to the deformation characterization of $T\ov\cM_{\be}$,
we find
\BE{cF_e}\begin{split}
c_1(\ov\cM_{\be}) &= -f_*\ev^*c_2(X) 
+\sum_{\underset{\be_1,\be_2\in H_+(X)}{\be_1+\be_2=\be}}
\!\!\!\!\!\!\!\!\!D_{\be_1,\be_2}\,,\\
2c_2(\ov\cM_{\be})-c_1^2(\ov\cM_{\be})
&=-f_*\big(\ev^*c_3(X)+\psi\,\ev_2^*c_2(X)+\psi^3\big)
+\frac{1}{2}\sum_{\underset{\be_1,\be_2\in H_+(X)}{\be_1+\be_2=\be}}
\!\!\!\!\!\!\!\!\!(\psi_1\!+\!\psi_2)D_{\be_1,\be_2} \,.
\end{split}\EE
The $4$-dimensional case of the first equation above appears 
in Section 1.2.4 of~\cite{KlP} and is also an immediate consequence of 
the $n\!=\!4$ analogue of~\e_ref{g1diff_e} below.
The second identity in~\e_ref{cF_e} is~\e_ref{g1diff_e} itself.

\subsubsection{The numbers \ref{d5n2p1psi0}}
\label{n2A_sssec}

\noindent
Suppose $\big((\bar\cC,x_e),(\cC^*,x_e)\big)$ is an element of~$\cZ'$.
Since the curve $\bar\cC\!\cup\!\cC^*$ passes through~$\mu$, 
$\bar\cC\!\cup\!\cC^*$ has at most two components.
We have three possibilities for~$\cZ'$.\\

\noindent
{\it Case $0$ ($\bar\cC\!=\!\cC^*$):}
Here $\be_1\!=\!\be_2$ and
$$\cZ'=\big\{\big((\cC^*,x_e),(\cC^*,x_e)\big)\!:
(\cC^*,x_e)\!\in\!\cM^*\big\}.$$
The normal bundle of $\cZ'$ in $\ov\cM\!\times\!\cM^*$ is isomorphic to
$T_{\be_1}\!\lra\!\cM^*$ and the differential
$$d\ev_{e,e}=d\ev_e\!: \N\lra \ev_{e,e}^*\N_{\De}$$
is injective over $\cM^*$. 
Thus, the contribution 
of $\cZ'$ to the homology intersection number is given~by
$$\blr{e(\N_{\De}/\N),\cZ'}=\blr{c_2(\N_{\be_1}),\cM^*}.$$
Using the second equation in~\e_ref{cN_e}, the first 
equation in \e_ref{cF_e}, and the fourth equation in \e_ref{n1psired_e},
we obtain the $\be_1\!=\!\be_2$ case of \e_ref{m2red_e1}.\\

\noindent
{\it Case $1A$ ($\bar\cC\!\subsetneq\!\cC^*$):}
Here $\be_1\!<\!\be_2$ and
$$\cZ'=\big\{\big((\bar\cC,x_e),(\bar\cC\!\vee\!\cC',x_e)\big)\!:
(\bar\cC,x_e)\!\in\!\ov\cM,\,
(\bar\cC\!\vee\!\cC',x_e)\!\in\!\cZ^*
\big\},$$
where $\cZ^*\!\subset\!\cM^*$ is the locus consisting of 2-component curves
with the marked point on the first component.
Thus, $\cZ'$ is the union of the first components of 
the finitely many $(\be_1,\be_2\!-\!\be_1)$-curves passing through
the constraint~$\mu$.
The normal bundle~$\N$ of $\cZ'$ in $\ov\cM\!\times\!\cM^*$ contains
the subbundle~$\pi_1^*T_{\be_1}$ and $\N/\pi_1^*T_{\be_1}$
is isomorphic to the normal bundle $\N\cZ^*$ of $\cZ^*$ in~$\cM^*$.
Since the differential
$$d\ev_{e,e}\!: \N\lra \ev_{e,e}^*\N_{\De}$$
is injective over $\cZ'$, the contribution of $\cZ'$ to 
the homology intersection number is given~by
$$\blr{e(\N_{\De}/\N),\cZ'}=\blr{c_1(\N_{\be_1})-c_1(\N\cZ^*),\cZ^*}.$$
Since the degrees of the restrictions of $\N_{\be_1}$ and $\N\cZ^*$ to 
each curve $\bar\cC$ are~$-2$ and~$-1$, respectively,
we obtain the $\be_1\!<\!\be_2$ case of~\e_ref{m2red_e1}.\\

\noindent
{\it Case $1B$ ($\bar\cC\!\supsetneq\!\cC^*$):}
Here $\be_1\!>\!\be_2$ and
$$\cZ'=\big\{\big((\cC'\vee\cC^*,x_e),(\cC^*,x_e)\big)\!:
(\cC'\!\vee\!\cC^*,x_e)\!\in\!\ov\cM,~
(\cC^*,x_e)\!\in\!\cZ^*\},$$
where $\cZ^*\!\subset\!\cM^*$ is the locus of curves meeting 
a $(\be_1\!-\!\be_2)$-curve.
Thus, $\cZ'$ is the union of the second components of 
the finitely many $(\be_1\!-\!\be_2,\be_2)$-curves whose second component
passes through the constraint~$\mu$.
The normal bundle~$\N$ of $\cZ'$ in $\ov\cM\!\times\!\cM^*$ contains
the subbundle~$\pi_1^*T_{\be_1}$ and $\N/\pi_1^*T_{\be_1}$
is isomorphic to the normal bundle $\N\cZ^*$ of $\cZ^*$ in~$\cM^*$.
The latter is trivial.
Since the differential
$$d\ev_{e,e}\!: \N\lra \ev_{e,e}^*\N_{\De}$$
is injective over $\cZ'$, the contribution of $\cZ'$ to 
the homology intersection number is given~by
$$\blr{e(\N_{\De}/\N),\cZ'}=\blr{c_1(\N_{\be_1})-c_1(\N\cZ^*),\cZ^*}.$$
The $\be_1\!>\!\be_2$ case of~\e_ref{m2red_e1} now follows
from the first equation in \e_ref{cN_e}.

\subsubsection{The numbers~\ref{d5n2p0psi0}}
\label{n2B_sssec}

\noindent
Suppose $\big((\bar\cC,x_e),(\cC^*,x_e)\big)$ is an element of~$\cZ'$.
The curve $\bar\cC\!\cup\!\cC^*$ then has one, two, or three components
and carries a marked point~$e$ lying on the divisor $\mu$.
The $6$ possibilities for the connected components of $\cZ'$
are indicated in Figure \ref{n5m2_fig}.\\

\noindent
{\it Case $0$ ($\bar\cC\!=\!\cC^*$):}
Here $\be_1\!=\!\be_2$ and
$\big((\bar\cC,x_e),(\cC^*,x_e)\big)$ is an element of
$$\ov\cS=\big\{\big((\bar\cC,x_e),(\bar\cC,x_e)\big)\!:
(\bar\cC,x_e)\!\in\!\ov\cM   \big\}\subset\cZ'.$$
The normal bundle of $\ov\cS$ in $\ov\cM\!\times\!\cM^*$ is isomorphic to
$T_{\be_2}\!\lra\!\ov\cM$, and 
the contribution of $\ov\cS$ to the homology intersection number is given by
$$\blr{e(\N_{\De}/\N),\ov\cS}=\blr{c_2(\N_{\be_2}),\ov\cM}.$$
Using  the second equation in \e_ref{cN_e} and 
the first equation in \e_ref{cF_e}, 
we obtain the $\be_1\!=\!\be_2$ case of the last term in 
\e_ref{rec_e2a}.\\

\noindent
{\it Case $1A$ ($|\cC^*|\!=\!2,~\bar\cC\!\subsetneq\!\cC^*$):}
Here $\be_1\!<\!\be_2$ and
$\big((\bar\cC,x_e),(\cC^*,x_e)\big)$ is an element~of
$$\ov\cS=\big\{\big((\bar\cC,x_e),(\bar\cC\!\vee\!\cC',x_e)\big)\!:
(\bar\cC,x_e)\!\in\!\ov\cZ,\,(\bar\cC\!\vee\!\cC',x_e)\!\in\!\cM^* 
\big\}
\subset\cZ',$$
where $\ov\cZ\!\subset\!\ov\cM$ is the locus consisting of curves meeting
a $(\be_2\!-\!\be_1)$-curve.
The normal bundle~$\N$ of $\ov\cS$ in $\ov\cM\!\times\!\cM^*$ contains
the subbundle~$\pi_2^*T_{\be_2}$ and $\N/\pi_2^*T_{\be_2}$
is isomorphic to the normal bundle~$\N\ov\cZ$ of~$\ov\cZ$ in~$\ov\cM$.
Since the differential
$$d\ev_{e,e}\!: \N\lra \ev_{e,e}^*\N_{\De}$$
is injective over $\ov\cS$, the contribution of $\ov\cS$ to 
the homology intersection number is given~by
$$\blr{e(\N_{\De}/\N),\ov\cS}=\blr{c_1(\N_{\be_2})-
\big(c_1(\N_{\be_2-\be_1})\!+\!\ti\psi_1\big),\ov\cZ},$$
where $\ti\psi_1$ is the untwisted $\psi$-class at the node of 
the  $(\be_2\!-\!\be_1)$-component of a curve in~$\ov\cZ$.

Using the first equations in \e_ref{cN_e} and in \e_ref{cF_e}
and the fourth equation in \e_ref{n1psired_e},
we obtain
 the $\be_1\!<\!\be_2$ case of the last term in \e_ref{rec_e2a} 
minus the last term in \e_ref{n2corr_e0}.
The latter arises from {\it Case~$2A$} below.\\

\noindent
{\it Case $1B$ ($|\bar\cC|\!=\!2,~\bar\cC\!\supsetneq\!\cC^*$):}
Here $\be_1\!>\!\be_2$  and
$\big((\bar\cC,x_e),(\cC^*,x_e)\big)$ is an element of
$$\ov\cS=\big\{\big((\cC'\vee\cC^*,x_e),(\cC^*,x_e)\big)\!:
(\cC'\vee\cC^*,x_e)\!\in\!\ov\cZ,\,(\cC^*,x_e)\!\in\!\cM^*
\big\}
\subset\cZ',$$
where $\ov\cZ\!\subset\!\ov\cM$ is 
the locus of $(\be_2,\be_1\!-\!\be_2)$-curves
with the marked point $e$ lying on the first 
component.
The normal bundle $\N$ of $\ov\cS$ in $\ov\cM\!\times\!\cM^*$ contains
the subbundle $\pi_2^*T_{\be_2}$, $\N/\pi_2^*T_{\be_2}$
is isomorphic to the normal bundle~$\N\ov\cZ$ of~$\ov\cZ$ in~$\ov\cM$, and
the contribution of $\cZ'$ to the homology intersection number is given~by
$$\blr{e(\N_{\De}/\N),\ov\cS}=\blr{c_1(\N_{\be_2})+
\big(\psi_1\!+\!\psi_2\big),\ov\cZ},$$
where $\psi_1$ and $\psi_2$ are the $\psi$-classes of the first and second 
components at the node of a curve in~$\ov\cZ$. 
We obtain
 the $\be_1\!>\!\be_2$ analogue of the {\it Case~$1A$} contribution 
in \e_ref{rec_e2a}.\\

\begin{figure}
\begin{pspicture}(-5,-6)(10,1.7)
\psset{unit=.4cm}
\psline[linewidth=.2](-2,-2)(-2,2)\pscircle*(-2,1){.25}
\rput(-2.7,-1.5){\sm{$\be_1$}}\rput(-1.4,1){\sm{$e$}}
\rput(-2,-3){\sm{$\be_2\!=\!\be_1$: Contr.~0}}
\psline[linewidth=.2](8,-2)(8,2)\psline(7,0)(10,3)
\pscircle*(8,0){.25}\rput(8.6,0){\sm{$e$}}
\rput(7.3,-1.5){\sm{$\be_1$}}\rput(11.5,3.2){\sm{$\be_2\!-\!\be_1$}}
\rput(8,-3){\sm{$\be_2\!>\!\be_1$: Contr.~1A}}
\psline(18,-2)(18,2)\psline[linewidth=.2](17,0)(20,3)
\pscircle*(19,2){.25}\rput(19.5,1.8){\sm{$e$}}
\rput(16.4,-1.5){\sm{$\be_1\!-\!\be_2$}}\rput(20.8,3.2){\sm{$\be_2$}}
\rput(18,-3){\sm{$\be_2\!<\!\be_1$: Contr.~1B}}
\psline[linewidth=.2](-2,-10)(-2,-6)
\psline(-3,-7.5)(0,-5)\pscircle*(-2,-8){.25}
\psline(-3,-8.5)(0,-11)
\rput(-2,-10.5){\sm{$\be_1$}}\rput(-1.4,-8){\sm{$e$}}
\rput(0.5,-4.8){\sm{$\be$}}\rput(0.5,-11.2){\sm{$\be'$}}
\rput(-2.3,-12.6){\sm{$\be\!+\!\be'\!=\!\be_2\!-\!\be_1$}}
\rput(-2,-14.2){\sm{$\be_2\!>\!\be_1$: Contr.~2A}}
\psline[linewidth=.2](8,-10)(8,-6)\psline(7,-7.5)(10,-5)
\pscircle*(8,-8){.25}
\psline(7,-8.5)(10,-11)
\rput(8,-10.5){\sm{$\be$}}\rput(8.6,-8){\sm{$e$}}
\rput(11.5,-4.8){\sm{$\be_2\!-\!\be$}}\rput(11.5,-11.2){\sm{$\be_1\!-\!\be$}}
\rput(8,-12.6){\sm{$\be\!<\!\be_1,\be_2$}}
\rput(8,-14.2){\sm{Contr.~2B}}
\psline[linewidth=.2](18,-10)(18,-6)\psline(17,-7.5)(20,-5)
\pscircle*(18,-8){.25}
\psline(17,-8.5)(20,-11)
\rput(18,-10.5){\sm{$\be_2$}}\rput(18.5,-8){\sm{$e$}}
\rput(20.5,-4.8){\sm{$\be$}}\rput(20.5,-11.2){\sm{$\be'$}}
\rput(17.7,-12.6){\sm{$\be\!+\!\be'\!=\!\be_1\!-\!\be_2$}}
\rput(18,-14.2){\sm{$\be_2\!<\!\be_1$: Contr.~2C}}
\end{pspicture}
\caption{Excess contributions for the meeting number $n_{\be_1\be_2}(\mu|;)$.
The labels refer to the cases described in 
Section \ref{n2B_sssec}.
The marked point $e$ corresponds to the (former) node and lies on the divisor~$\mu$.
The thicker lines indicate the multiple component.
The space of curves in the first diagram in the top row is $2$-dimensional.
The other two spaces in the first row are $1$-dimensional.
All spaces in the bottom row are $0$-dimensional.}
\label{n5m2_fig}
\end{figure}
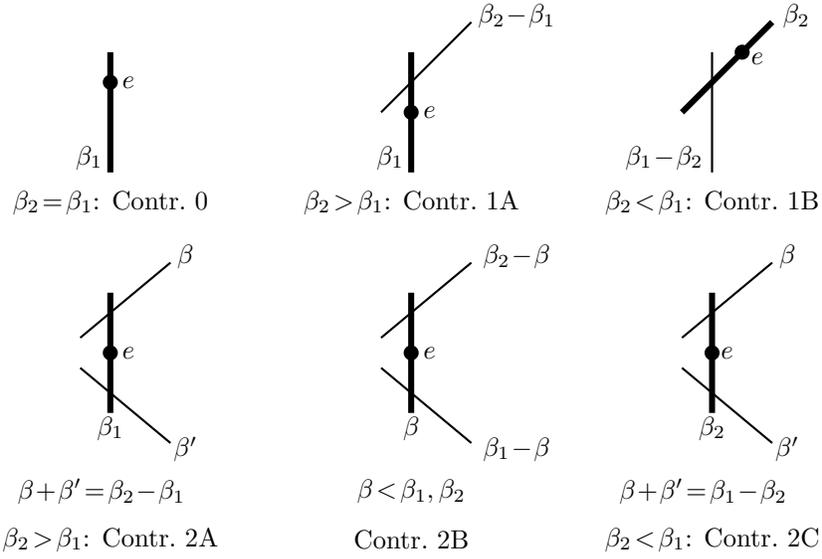

\noindent
{\it Case $2A$ ($|\cC^*|\!=\!3,~\bar\cC\!\subsetneq\!\cC^*$):}
Here $\be_1\!<\!\be_2$.
If $|\bar\cC|\!=\!2$, $\big((\bar\cC,x_e),(\cC^*,x_e)\big)$ is an element
of the space $\bar\cS$ in {\it Case~$1A$} above. 
This is also the case if $|\bar\cC|\!=\!1$ and 
the curve $\cC^*\!-\!\bar\cC$ is connected.
In the remaining case, $\bar\cC$ is the middle component of the 3-component
curve~$\cC^*$ and carries the marked point~$e$, which lies on the
divisor $\mu$.
Each such pair $\big((\bar\cC,x_e),(\cC^*,x_e)\big)$ is a regular element 
of $\cZ$ and therefore contributes~$1$ to the homology intersection.
The contribution of such pairs is accounted for by 
the last term in \e_ref{n2corr_e0}.\\

\noindent
{\it Case $2B$ ($|\bar\cC|\!=\!|\cC^*|\!=\!2,~\bar\cC\!\neq\!\cC^*$):}
Here the curve $\bar\cC \cup \cC^*$ consists of three components,
with the middle component meeting the hyperplane 
$\mu$ at the marked point $e$.
Such pairs $\big((\bar\cC,x_e),(\cC^*,x_e)\big)$ are regular elements
of~$\cZ$, and their contribution is accounted for by 
the middle term on the right side of \e_ref{rec_e2a}.\\

\noindent
{\it Case $2C$ ($|\cC^*|\!=\!3,~\bar\cC\!\supsetneq\!\cC^*$):}
The analysis
is the same as  {\it Case~$2A$} with $\be_1$ and~$\be_2$ interchanged.

\subsubsection{The numbers \ref{d5n3}}
\label{n3_sssec}

\noindent
If $\big((\bar\cC,x_e),(\cC^*,x_e)\big)$ is an element of~$\cZ'$,
$\bar\cC$ consists of two sets of components, $\bar\cC_1$
and~$\bar\cC_2$, with the second component carrying the marked point~$e$.
Either $\bar\cC_1$ or~$\bar\cC_2$ may consist of two components,
while the other curve must consist of one component.
The total number of components in $\bar\cC \cup \cC^*$ is either 
two or three.
The $12$ possibilities for the connected components of $\cZ'$
are indicated in Figure~\ref{n5m3_fig}.\\

\noindent
{\it Case $0$ ($|\bar\cC\!\cup\!\cC^*|\!=\!2,~\bar\cC_2\!\subset\!\cC^*$):}
If $\cC^*\!=\!\bar\cC_2$, then $\be_2\!=\!\be_3$ and
$\big((\bar\cC,x_e),(\cC^*,x_e)\big)$ is an element~of
$$\ov\cS=\big\{\big((\bar\cC_1\vee\cC^*,x_e),(\cC^*,x_e)\big)\!:
(\bar\cC_1\vee\cC^*,x_e)\!\in\!\ov\cM\big\}  \subset\cZ'.$$
Similarly to {\it Case~$0$} 
in Sections~\ref{n2A_sssec} and~\ref{n2B_sssec},  
the normal bundle of $\ov\cS$ in $\ov\cM\!\times\!\cM^*$ is  isomorphic to
$T_{\be_2}\!\lra\!\ov\cM$, and 
the contribution of $\ov\cS$ to the homological intersection number is given~by
$$\blr{e(\N_{\De}/\N),\ov\cS}=\blr{c_2(\N_{\be_2}),\ov\cM}.$$
Using the second equation in \e_ref{cN_e}, 
the first equation in \e_ref{cF_e}, and 
the fourth equation in \e_ref{n1psired_e},
we obtain
 the $\be_3\!=\!\be_2$ case of the term $\tC_{\be_1\be_2\be_3}^{(2)}$
in \e_ref{rec_e2b}.\\

\noindent
If $\cC^*\!=\!\bar\cC_1\cup\bar\cC_2$, then $\be_1\!+\be_2\!=\!\be_3$ and
$\big((\bar\cC,x_e),(\cC^*,x_e)\big)$ is an element~of
$$\ov\cS=\big\{\big((\cC^*,x_e),(\cC^*,x_e)\big)\!:
(\cC^*,x_e)\!\in\!\ov\cM\big\}
\subset\cZ'.$$
The normal bundle of $\ov\cS$ in $\ov\cM\!\times\!\cM^*$ is  isomorphic to
$T_{\be_1+\be_2}\!\lra\!\ov\cM$, and 
the contribution of $\ov\cS$ to the homological 
intersection number is given by
$$\blr{e(\N_{\De}/\N),\ov\cS}=\blr{c_2(\N_{\be_1+\be_2}),\ov\cM}.$$
Using  the second equation in \e_ref{cN_e},
the first equation in \e_ref{cF_e}, and 
the fourth equation in~\e_ref{n1psired_e},
we obtain the $\be_3\!=\!\be_1\!+\be_2$ case of the term 
$\tC_{\be_1\be_2\be_3}^{(12)}$ in \e_ref{rec_e2b}.\\

\noindent
{\it Case $0'$ 
($|\bar\cC\!\cup\!\cC^*|\!=\!2,~\bar\cC_2\!\not\subset\!\cC^*$):}
Here  $\be_1\!=\!\be_3$ and
$\big((\bar\cC,x_e),(\cC^*,x_e)\big)$ is an element~of
$$\ov\cS=\big\{\big((\cC^*\vee\bar\cC_2,x_e),(\bar\cC,x_e)\big)\!:
(\cC^*\vee\bar\cC_2,x_e)\!\in\!\ov\cZ\big\}  \subset\cZ',$$
where $\ov\cZ\!\subset\!\ov\cM$ consists of the pairs of $1$-marked curves
with the marked point at the node of the two curves.
The normal bundle $\N$ of $\ov\cS$ in $\ov\cM\!\times\!\cM^*$ contains 
$T_{\be_1}$ as a subbundle, and $\N/T_{\be_1}$ is isomorphic to the normal bundle 
of $\ov\cZ$ in~$\ov\cM$. The latter is the universal tangent line bundle
at the marked point.
Since the homomorphism  
$$d\ev_{e,e}\!: \N\lra \ev_{e,e}^*\N_{\De}$$
is injective over $\cZ'$, the contribution of $\cZ'$ to 
the homological intersection number is given by
$$\blr{e(\N_{\De}/\N),\ov\cS}=\blr{c_1(\N_{\be_2})+\psi_2,\ov\cZ}.$$
Using  the first equations in \e_ref{cN_e} and \e_ref{cF_e} and 
the fourth equation in \e_ref{n1psired_e},
we obtain the $\be_3\!=\!\be_1$ case of the term $\tC_{\be_1\be_2\be_3}^{(1)}$
in \e_ref{rec_e2b}.\\

\begin{figure}
\begin{pspicture}(-6,-11)(10,1.7)
\psset{unit=.4cm}
\psline[linewidth=.2](-9,-2)(-9,2)\psline(-10,-.5)(-7,-3)\pscircle*(-9,.5){.25}
\rput(-9,-2.5){\sm{$\be_2$}}\rput(-8.4,.5){\sm{$e$}}
\rput(-6.4,-3.2){\sm{$\be_1$}}
\rput(-8.5,-5){\sm{$\be_3\!=\!\be_2$: Contr.~$0$}}
\psline[linewidth=.2](0,-2)(0,2)\psline(-1,.5)(2,3)\psline(-1,-.5)(2,-3)
\pscircle*(0,0){.25}\rput(.5,0){\sm{$e$}}
\pscircle*(0,1.33){.2}\rput(-.5,1.8){\sm{$e'$}}
\rput(0,-2.5){\sm{$\be_2$}}
\rput(2.6,-3.2){\sm{$\be_1$}}\rput(3.5,3.2){\sm{$\be_3\!-\!\be_2$}}
\rput(0.5,-5){\sm{$\be_3\!>\!\be_2$: Contr.~$1A$}}
\psline[linewidth=.2](9,-2)(9,2)\psline(8,.5)(11,3)\pscircle*(9,1.33){.25}
\psline(8,-.5)(11,-3)
\rput(9,-2.6){\sm{$\be_1$}}\rput(9.5,1.2){\sm{$e$}}
\rput(11.7,3.2){\sm{$\be_2$}}\rput(12.3,-3.2){\sm{$\be_3\!-\!\be_1$}}
\rput(10,-5){\sm{$\be_3\!>\!\be_1$: Contr.~$1A'$}}
\psline[linewidth=.2](18,-2)(18,2)\psline(17,.5)(20,3)
\psline[linewidth=.2](17,-0.5)(20,-3)
\pscircle*(18,0){.25}\rput(18.5,0){\sm{$e$}}
\pscircle*(18,1.33){.2}\rput(17.5,1.8){\sm{$e'$}}
\rput(17.5,-2.6){\sm{$\be_3\!-\!\be_1$}}
\rput(22.4,3.2){\sm{$\be_1\!+\!\be_2\!-\!\be_3$}}\rput(20.6,-3.2){\sm{$\be_1$}}
\rput(19.5,-5){\sm{$\be_1\!<\!\be_3\!<\!\be_1\!+\!\be_2$:}}
\rput(18,-6){\sm{Contr.~$1B$}}
\psline[linewidth=.2](-9,-12)(-9,-8)\psline[linewidth=.2](-10,-10.5)(-7,-13)
\pscircle*(-9,-9.5){.25}
\rput(-9,-12.5){\sm{$\be_2$}}\rput(-8.4,-9.5){\sm{$e$}}
\rput(-6.4,-13.2){\sm{$\be_1$}}
\rput(-8.5,-15){\sm{$\be_3\!=\!\be_1\!+\!\be_2$: Contr.~$0$}}
\psline[linewidth=.2](0,-12)(0,-8)\psline(-1,-9.5)(2,-7)
\psline[linewidth=.2](-1,-10.5)(2,-13)
\pscircle*(0,-10){.25}\rput(.6,-10){\sm{$e$}}
\pscircle*(0,-8.67){.2}\rput(-.5,-8.2){\sm{$e'$}}
\rput(0,-12.5){\sm{$\be_2$}}
\rput(2.6,-13.2){\sm{$\be_1$}}\rput(4.4,-6.8){\sm{$\be_3\!-\!\be_1\!-\!\be_2$}}
\rput(1,-15){\sm{$\be_3\!>\!\be_1\!+\!\be_2$: Contr.~$1A$}}
\psline[linewidth=.2](9,-12)(9,-8)\psline(8,-9.5)(11,-7)\psline(8,-10.5)(11,-13)
\pscircle*(9,-10){.25}\rput(9.6,-10){\sm{$e$}}
\pscircle*(9,-8.77){.2}\rput(8.5,-8.2){\sm{$e'$}}
\rput(9,-12.5){\sm{$\be_3$}}
\rput(12.5,-6.8){\sm{$\be_2\!-\!\be_3$}}\rput(11.6,-13.2){\sm{$\be_1$}}
\rput(10,-15){\sm{$\be_3\!<\!\be_2$: Contr.~$1B$}}
\psline[linewidth=.2](18,-12)(18,-8)\psline[linewidth=.2](17,-9.5)(20,-7)
\psline(17,-10.5)(20,-13)
\pscircle*(19,-7.83){.25}\rput(19.4,-8.3){\sm{$e$}}
\pscircle*(18,-11.3){.2}\rput(18.5,-10.8){\sm{$e'$}}
\rput(17.5,-12.6){\sm{$\be_3\!-\!\be_2$}}
\rput(22.2,-13.2){\sm{$\be_1\!+\!\be_2\!-\!\be_3$}}\rput(20.6,-7.2){\sm{$\be_2$}}
\rput(19.5,-15){\sm{$\be_2\!<\!\be_3\!<\!\be_1\!+\!\be_2$:}}
\rput(19,-16){\sm{Contr.~$1B$}}
\psline(-9,-23)(-9,-19)\psline[linewidth=.2](-10,-21.5)(-7,-24)
\pscircle*(-9,-22.33){.25}
\rput(-9,-23.5){\sm{$\be_2$}}\rput(-8.4,-22.1){\sm{$e$}}
\rput(-6.4,-24.2){\sm{$\be_1$}}
\rput(-8.5,-26){\sm{$\be_3\!=\!\be_1$: Contr.~$0'$}}
\psline[linewidth=.2](0,-23)(0,-19)\psline[linewidth=.2](-1,-20.5)(2,-18)
\psline(-1,-21.5)(2,-24)
\pscircle*(1,-18.83){.25}\rput(1.2,-19.3){\sm{$e$}}
\pscircle*(0,-22.3){.2}\rput(.5,-21.8){\sm{$e'$}}
\rput(0,-23.5){\sm{$\be_1$}}
\rput(2.6,-17.8){\sm{$\be_2$}}\rput(4.2,-24.2){\sm{$\be_3\!-\!\be_1\!-\!\be_2$}}
\rput(1.5,-26){\sm{$\be_3\!>\!\be_1\!+\!\be_2$: Contr.~$1A$}}
\psline(9,-23)(9,-19)\psline[linewidth=.2](8,-20.5)(11,-18)
\psline(8,-21.5)(11,-24)
\pscircle*(10,-18.83){.25}\rput(10.4,-19.2){\sm{$e$}}
\pscircle*(9,-22.33){.2}\rput(9.5,-21.8){\sm{$e'$}}
\rput(8.5,-23.6){\sm{$\be_2\!-\!\be_3$}}
\rput(11.7,-17.9){\sm{$\be_3$}}\rput(11.6,-24.2){\sm{$\be_1$}}
\rput(10.5,-26){\sm{$\be_3\!<\!\be_2$: Contr.~$1B$}}
\psline[linewidth=.2](18,-23)(18,-19)
\psline(17,-20.5)(20,-18)\pscircle*(18,-19.67){.25}\psline(17,-21.5)(20,-24)
\rput(18,-23.6){\sm{$\be_3$}}\rput(18.6,-19.8){\sm{$e$}}
\rput(20.7,-17.8){\sm{$\be_2$}}\rput(21.3,-24.2){\sm{$\be_1\!-\!\be_3$}}
\rput(19.5,-26){\sm{$\be_3\!<\!\be_1$: Contr.~$1B'$}}
\end{pspicture}
\caption{Excess contributions for the meeting number $m_{\be_1\be_2\be_3}$.
The labels  refer to the cases described in 
Section \ref{n3_sssec}.
The marked point $e$ corresponds to the (former) node joining
the $\be_2$ and $\be_3$ curves.
For the curves of types~1A and~1B,
$e'$ indicates the new node on the (leftover) $(\be_1,\be_2)$-curve.
The thicker lines represent the multiple component(s).
The excess loci corresponding to {\it Contr.~$0$} are $2$-dimensional.
The loci corresponding to {\it Contr.~$1A'$} and $1B'$ are $0$-dimensional.
The remaining loci are $1$-dimensional.}
\label{n5m3_fig}
\end{figure}
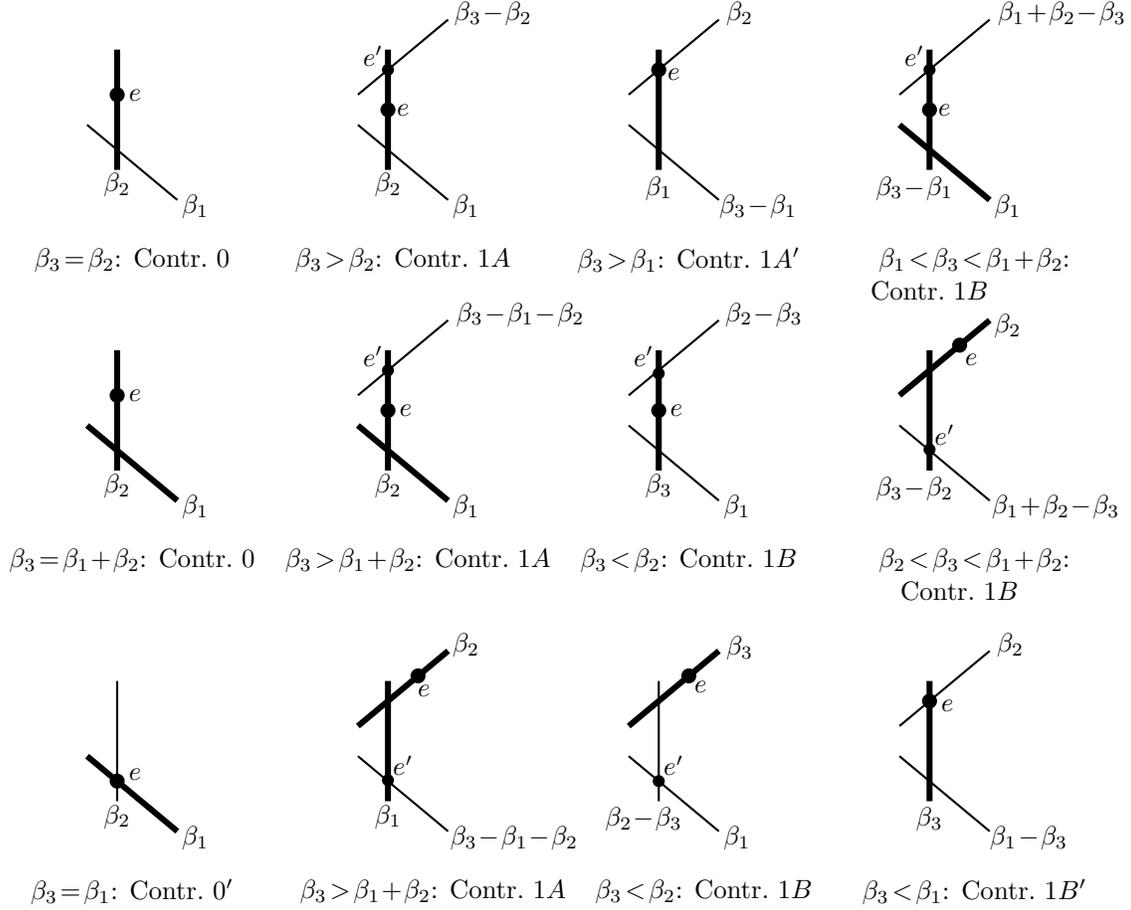

\noindent
{\it Case $1A$ ($|\bar\cC \cup \cC^*|\!=\!3,~\cC^*\!\not\subset\!\bar\cC,~
\bar\cC_2\!\subset\!\cC^*$):}
If $\bar\cC_1\!\not\subset\!\cC^*$, then $\be_2\!<\!\be_3$ and
$\big((\bar\cC,x_e),(\cC^*,x_e)\big)$ is an element~of
$$\ov\cS=\big\{\big((\bar\cC_1\vee\bar\cC_2,x_e),(\bar\cC_2\vee\cC',x_e)\big)\!:
(\bar\cC_1\vee\bar\cC_2,x_e)\!\in\!\ov\cZ\big\}  \subset\cZ',$$
where $\ov\cZ\!\subset\!\ov\cM$ consists of the pairs of $(\be_1,\be_2)$-curves
such that the second component meets a $(\be_3\!-\!\be_2)$-curve.
We see
$\ov\cS$ is the union of the middle components of 
$(\be_1,\be_2,\be_3\!-\!\be_2)$-curves in~$X$,
with each curve contributing $-1$ to the homological intersection number.
The contribution accounts for the $\be_3\!>\!\be_2$ case of 
the term $\tC_{\be_1\be_2\be_3}^{(2)}$ in~\e_ref{rec_e2b}.\\

\noindent
If $\bar\cC_1\!\subset\!\cC^*$, then $\be_1\!+\!\be_2\!<\!\be_3$ and
$\big((\bar\cC,x_e),(\cC^*,x_e)\big)$ is an element~of
$$\ov\cS=\big\{\big((\bar\cC_1\vee\bar\cC_2,x_e),
(\bar\cC_1\vee\bar\cC_2\vee\cC',x_e)\big)\!:
(\bar\cC_1\vee\bar\cC_2,x_e)\!\in\!\ov\cZ\big\}  \subset\cZ',$$
where $\ov\cZ\!\subset\!\ov\cM$ consists of the pairs $(\be_1,\be_2)$-curves
meeting a $(\be_3\!-\!\be_1\!-\!\be_2)$-curve with 
the $\be_2$-component carrying the marked point~$e$.
Here, $\bar\cS$ is the union of the last components of
the $(\be_3\!-\!\be_1\!-\!\be_2,\be_1,\be_2)$-curves and 
the middle components of $(\be_1,\be_2,\be_3\!-\!\be_1\!-\!\be_2)$-curves.
By reasoning analogous to {\it Case~$1A$} in Section~\ref{n2B_sssec},
each of the former contributes $-1$ to the homological intersection number,
while each of the latter contributes $0$.
We obtain the $\be_3\!>\!\be_1\!+\!\be_2$ case of 
the term $\tC_{\be_1\be_2\be_3}^{(12)}$ in~\e_ref{rec_e2b}.\\

\noindent
{\it Case $1A'$ ($|\bar\cC\!\cup\!\cC^*|\!=\!3,~\cC^*\!\not\subset\!\bar\cC,~
\bar\cC_2\!\not\subset\!\cC^*$):}
Here $\be_3\!>\!\be_1$ and $(\bar\cC \cup \cC^*,x_{e})$
is a $(\be_3\!-\!\be_1,\be_1,\be_2)$-curve with the marked point lying 
on the node joining the last two components.
Each such pair $\big((\bar\cC,x_e),(\cC^*,x_e)\big)$ is a regular element of~$\cZ$,
contributing~$1$ to the homology intersection number.
We obtain  the $\be_3\!>\!\be_1$ case of the term 
$\tC_{\be_1\be_2\be_3}^{(1)}$ in~\e_ref{rec_e2b}.\\

\noindent
{\it Case $1B$ ($|\bar\cC\!\cup\!\cC^*|\!=\!3,~\cC^*\!\subset\!\bar\cC,~
\cC^*\!\not\subset\!\bar\cC_1$):}
If $\cC^*\!=\!\bar\cC_2$ or $\cC^*\!=\!\bar\cC$, $\big((\bar\cC,x_e),(\cC^*,x_e)\big)$ 
is an element of 
one of the spaces $\bar\cS$ defined in {\it Case~$0$} above.
Hence, we can assume that $\cC^*\!\neq\!\bar\cC_2,\bar\cC$.
If $\cC^*\!\subset\!\bar\cC_2$, then $\be_2\!>\!\be_3$ and
$\big((\bar\cC,x_e),(\cC^*,x_e)\big)$ is an element~of
$$\ov\cS=\big\{\big((\bar\cC_1\vee\cC^*\vee\cC',x_e),(\cC^*,x_e)\big)\!:
(\bar\cC_1\vee\cC^*\vee\cC',x_e)\!\in\!\ov\cZ\big\}  \subset\cZ',$$
where $\ov\cZ\!\subset\!\ov\cM$ is the locus of the pairs $(\be_1,\be_2)$-curves
with the second component broken into two.
As in {\it Case~1B} of Section \ref{n2B_sssec},
$\ov\cS$ is the union of the middle components of 
$(\be_1,\be_3,\be_2\!-\!\be_3)$-curves and 
the last components of $(\be_1,\be_2\!-\!\be_3,\be_3)$,
with each curve contributing~$-1$ to the homological intersection number.
The contribution accounts for the $\be_3\!<\!\be_2$ case of 
the term $\tC_{\be_1\be_2\be_3}^{(2)}$ in~\e_ref{rec_e2b}.\\

\noindent
If $\cC^*\!\not\subset\!\bar\cC_2$, then $\be_1\!+\!\be_2\!>\!\be_3$ and
$\big((\bar\cC,x_e),(\cC^*,x_e)\big)$ is an element of
$$\ov\cS=\big\{\big((\bar\cC_1\vee\bar\cC_2\vee\cC',x_e),
(\bar\cC_1\vee\bar\cC_2,x_e)\big)\!:
(\bar\cC_1\vee\bar\cC_2\vee\cC',x_e)\!\in\!\ov\cZ\big\}  \subset\cZ',$$
where $\ov\cZ\!\subset\!\ov\cM$ is the locus of the 
pairs $(\be_1,\be_2)$-curves
with one of the components broken into two.
Here, $\ov\cS$ is the union of the middle components of 
$(\be_1,\be_3\!-\!\be_1,\be_1\!+\!\be_2\!-\!\be_3)$-curves,
if $\be_3\!>\!\be_1$,  and the last components of 
$(\be_1\!+\!\be_2\!-\!\be_3,\be_3\!-\!\be_2,\be_2)$-curves, if $\be_3\!>\!\be_2$.
Each of the latter curves contributes $-1$  to the homological
 intersection number,
while each of the former contributes $0$.
We obtain the $\be_3\!<\!\be_2$ case of 
the term $\tC_{\be_1\be_2\be_3}^{(12)}$ in~\e_ref{rec_e2b}.\\

\noindent
{\it Case $1B'$ ($|\bar\cC\!\cup\!\cC^*|\!=\!3,~\cC^*\!\subset\!\bar\cC_1$):}
If $\cC^*\!=\!\bar\cC_1$, then $\big((\bar\cC,x_e),(\cC^*,x_e)\big)$ 
is an element~of  the space~$\bar\cS$ defined in {\it Case~$0'$} above.
Hence, we can assume that $\cC^*\!\neq\!\bar\cC_1$.
Then, $\be_1\!>\!\be_3$ and
$\big((\bar\cC,x_e),(\cC^*,x_e)\big)$ is an element~of
$$\ov\cS=\big\{\big((\cC'\vee\cC^*\vee\bar\cC_2,x_e),(\cC^*,x_e)\big)\!:
(\cC'\vee\cC^*\vee\bar\cC_2,x_e)\!\in\!\ov\cZ\big\}  \subset\cZ',$$
where $\ov\cZ\!\subset\!\ov\cM$ is the locus of the pairs of 
$1$-marked $(\be_1,\be_2)$-curves represented by 
a $(\be_1\!-\!\be_3,\be_3,\be_2)$-curve in $X$ with the marked point
on the node of the last two components.
Each such pair $\big((\bar\cC,x_e),(\cC^*,x_e)\big)$ 
is a regular element of $\cZ$,
contributing $1$ to the homological intersection number and 
accounting for the $\be_3\!<\!\be_1$ case of 
the term $\tC_{\be_1\be_2\be_3}^{(1)}$ in~\e_ref{rec_e2b}.

\section{Genus 1 counts}
\label{g1nums_sec}

\subsection{Overview}
\label{g1dfn_subs}

\noindent
For each $\be\!\in\!H_+(X)$, $N_{1,\be}$ is the number of 
automorphism-weighted stable  
$C^{\i}$-maps  
$$u\!:\Si\!\lra\!X$$ from prestable curve of genus 1 to
$X$ of degree $\be$
solving 
 a perturbed Cauchy-Riemann equation,
\BE{CR_e} \dbar u+\nu(u)=0, \EE
for a small generic multi-valued perturbation $\nu$,
see Section 1.3 of \cite{g1comp2} for more details.
If $X$ is an ideal Calabi-Yau 
 $n$-fold, $\ov\M_1(X,\be)$ decomposes into strata
$\cZ_{\T}$ which each have well-defined
contribution to $N_{1,\be}$
in following sense:
\begin{quotation}\begin{it}
\noindent
For every stratum $\cZ_{\T}$, there exist $\tC_{\T}(\be)\!\in\!\Q$,  
$\ep_{\nu}\!\in\!\R^+$, and a compact subset $K_{\nu}$ of $\cZ_{\T}$ 
with the following property.
For every compact subset $K$ of $\cZ_{\T}$ and an open neighborhood~$U$ 
of~$K$ in the space of stable $C^\infty$-maps,
there exist an open neighborhood $U_{\nu}(K)$ of $K$ 
and $\ep_{\nu}(U)\!\in\!(0,\ep_{\nu})$, 
respectively,\footnote{$U_{\nu}(K)$ depends on $K$, 
while $\ep_{\nu}(U)$ depends on $U$.} such that
$$^{\pm}\big|\{\dbar\!+\!t\nu\}^{-1}(0)\!\cap\!U\big|
=\tC_{\T}(\be)
\quad\hbox{if}~~
t\!\in\!(0,\ep_{\nu}(U)),~K_{\nu}\!\subset\! 
K\!\subset\! U\!\subset\! U_{\nu}(K).$$
\end{it}\end{quotation}
While there are many different strata, it turns out that
$\tC_{\T}(\be)\!\neq\!0$ only for strata of the three simplest types.

If $X$ is an ideal Calabi-Yau $n$-fold, 
there are finitely many genus 1 curves in each homology class of~$X$. 
Furthermore,  every genus 1 curve $\cC$ in $X$ is embedded and super-rigid:
if $\N$ is the normal bundle of $\cC$ and 
$$u\!:\Si\!\lra\!\cC$$ is an unramified cover, 
then $H^0(\Si,u^*\N)\!=\!0$.
Hence, $H^1(\Si,u^*\N)\!=\!0$ and for every $d\!\in\!\Z^+$
$$\cZ_{(1,\be/d)} = \bigcup_{[\cC]=\be/d}\ov\M_1(\cC,d)$$
is a finite set of isolated regular points of $\ov\M_1(X,\be)$.
Each such point $u$ contributes $|1/\Aut(u)|$ to $N_{1,\be}$.
If $n_{1,\be}$ is the number of genus 1 curves in the homology class $\be$,
then
\BE{BPScontr_e} \tC_{(1,\be/d)}(\be)=\frac{\si(d)}{d}n_{1,\be/d}, \EE
where $\si(d)$ is the number of degree $d$ unbranched covers of
a genus 1 curve by connected genus 1 curves.
The integral number $n_{1,\be/d}$ is zero unless $d|\be$, or equivalently,
$\be/d$ is an integral homology class.

The remaining elements $u\!:\Si\!\lra\!X$ of $\ov\M_1(X,\be)$
are maps to genus 0 curves in $X$.
They split into strata $\cZ_{\T}$ indexed by combinatorial data
described in Section \ref{g1summ_subs}.
We will call a stratum $\cZ_{\T}$ {\em basic} if either of the
following conditions holds:
\begin{enumerate}[label=(B\arabic*)]
\item\label{basiceff_item} 
the domain $\Si$ of every element $[\Si,u]$ of $\cZ_{\T}$ is 
a nonsingular genus 1 curve, or
\item\label{basicghost_item} the domain $\Si$ of every element $[\Si,u]$ 
of $\cZ_{\T}$ is a union of a nonsingular genus 1 curve $\Si_P$ and a $\bP^1$
and $u$ is constant on $\Si_P$.
\end{enumerate}
In both cases, the restriction of $u$ to the non-contracted component
must be a $d\!:\!1$ cover of a curve in the homology class $\be/d$,
for some $d\!\in\!\Z^+$.
We will write $\T_{\eff}(\be/d,d)$ and $\T_{\gh}(\be/d,d)$
for the corresponding types of strata 
\ref{basiceff_item} and \ref{basicghost_item},
with {\em eff} and {\em gh} standing for {\it effective} and {\it ghost}
(principal component).

\begin{thm}
\label{g1thm}
Suppose $X$ is an ideal Calabi-Yau $5$-fold.
\begin{enumerate}
\item[(i)] If $\cZ_{\T}$ is a stratum of $\ov\M_1(X,\be)$ consisting of maps
to rational curves in $X$ and is not basic, $\tC_{\T}(\be)\!=\!0$.
\item[(ii)] 
For $\be\!\in\!H_+(X)$ and $d\!\in\!\Z^+$, 
\BE{g1thm_e}
\tC_{\T_{\eff}(\be,d)}(d\be)=\frac{d-1}{d^2}\tC_{\T_{\gh}(\be,1)}(\be),
\qquad 
\tC_{\T_{\gh}(\be,d)}(d\be)=\frac{1}{d^2}\tC_{\T_{\gh}(\be,1)}(\be).\EE
\end{enumerate}
\end{thm}

In Section \ref{g1ghost_subs}, we will prove
\BE{g1deg1contr_e}
\tC_{\T_{\gh}(\be,1)}(\be)=
\frac{1}{24}\int_{\ov\cM_{\be}}
\big(2c_2(\ov\cM_{\be})-c_1^2(\ov\cM_{\be})\big).\EE
On the other hand, the space $\ov\cM_\beta=\ov\M_0^*(X,\be)$
consists of regular maps to $X$. 
Thus, the contribution to $N_{1,\be}$ is given by the right side 
of equation~(2.15) in \cite{g1diff}:
$$\tC_{\T_{\gh}(\be,1)}(\be)
=\frac{1}{24}\bigg(
-n_{\be}\Big(\frac{c(X)}{1-\psi}\Big)
+\frac{1}{2}\sum_{\underset{\be_1,\be_2\in H_+(X)}{\be_1+\be_2=\be}}
\!\!\!\!\!\!n_{\be_1\be_2}\big(\psi_1\!+\!\psi_2|;)\bigg).$$
Comparing the above identity with \e_ref{g1deg1contr_e}, we find that 
\BE{g1diff_e} 
\int_{\ov\cM_{\be}}
\big(2c_2(\ov\cM_{\be})-c_1^2(\ov\cM_{\be})\big)
=-n_{\be}\Big(\frac{c(X)}{1-\psi}\Big)
+\frac{1}{2}\sum_{\underset{\be_1,\be_2\in H_+(X)}{\be_1+\be_2=\be}}
\!\!\!\!\!\!n_{\be_1\be_2}(\psi_1\!+\!\psi_2|;).\EE
We calculate the left side in terms of the Gromov-Witten invariants
of $X$ by expanding the right side via the equations of Section~\ref{g0nums_sec}.

Our proof of Theorem \ref{g1thm} applies also 
in dimensions $3$ and $4$. 
In particular, the result provides a direct explanation of the $1/d$-scaling
in the latter case discovered by other means in Section~$2$ of~\cite{KlP}.
Many aspects of the proof are applicable in dimensions $6$ and higher as well.

\subsection{Preliminaries}
\label{g1summ_subs}

Let $X$ be an ideal Calabi-Yau $5$-fold.
The strata of $\ov\M_1(X,\be)$ consisting of maps to rational curves 
can be described by {\em decorated graphs} 
$$\T=\big(\Ver,\Edg,\d,\un\be,\ka,i^*\big),$$
where
\begin{enumerate}[label=(D\arabic*)]
\item $\Ga\!=\!(\Ver,\Edg)$ is a connected graph 
containing 
either exactly one loop or a distinguished vertex, but not both,
\item $\un\be\!=\!(\be_i)_{i\in[m]}$ is an $m$-tuple of 
elements of~$H_+(X)$, with $m\!\in\!\{1,2,3\}$,
\item $\d\!:\Ver\!\lra\!\Z^{\geq 0}$ is a map,
$\ka\!:\d^{-1}(\Z^+)\!\lra\![m]$ is a surjective map,
\item  $i^*\!\in\!\{\star\}\cup[m]$.
\end{enumerate}
The irreducible components and the nodes of the domain $\Si$ of every element 
$[\Si,u]$ of $\cZ_{\T}$ correspond to the sets $\Ver$ and $\Edg$ respectively.
If $v\!\in\!\Ver$ is not the distinguished vertex of $\Ga$,
the corresponding component $\Si_v$ of $\Si$ is a $\P$.
Otherwise, $\Si_v$ is nonsingular of genus 1.
If $v\!\in\!\Ver$, the restriction of $u$ to $\Si_v$ is 
constant if $\d(v)\!=\!0$.
If $\d(v)\!\neq\!0$, $u|_{\Si_v}$ is a $\d(v)\!:\!1$ cover of
the component $\cC_{\ka(v)}$ of $\cC$.
If $\d$ does not vanish identically of the loop
in the graph $(\Ver,\Edg)$ or on the distinguished
vertex, $i^*$ is set to $\star$. If $\d$ vanishes
identically on the loop or on the  distinguished vertex,
the corresponding components of $\Si$ are mapped by $u$ to
a point on the $i^*$-component of $\cC$.
Since $u$ is continuous, $\cZ_{\T}\!=\!\eset$ unless 
$\ka$ satisfies certain combinatorial conditions.{\footnote{The 
strata $\cZ_{\T}$ as defined above intersect if $m\!\ge\!2$
and $\d$ vanishes on the loop or the distinguished vertex of $(\Ver,\Edg)$.
The issue can be easily addressed by allowing $i^*$ to take values 
in $\{\star\} \cup [m]\!\cup\!\{(1,2),(2,3)\}$.
However, equation~\e_ref{CR_e2} will be shown to
have no solutions near $\cZ_{\T}$ for a good choice of $\nu$ if $m\!\ge\!2$,
so further discussion is not needed.}}

Given a generic deformation of $\nu$ of the $\dbar$-operator as 
in~\e_ref{CR_e} and sufficiently small $t\!\in\!\R^+$, 
we will determine the number of solutions $[\Si,u]$ of
\BE{CR_e2} \dbar u+t\nu(u)=0, \EE
with $u$ close to the stratum $\cZ_{\T}$.
The assumption that $\nu$ is generic implies
that all solutions of \e_ref{CR_e2} are maps from nonsingular genus 1
 curves.
The arguments follow \cite{g2n2and3,g1comp}.
In particular, the gluing construction for $\cZ_{\T}$ will be performed
on a family of representatives $(\Si,u)$ for the elements $[\Si,u]$
in $\cZ_{\T}$, see Section 2.2 of \cite{g1comp}.
Our treatment here is less explicit
in order to 
streamline the discussion.

For the rest of Section \ref{g1nums_sec}, we fix a 
decorated graph $\T$ as above.
We define
$$|\un\be|=\sum_{i=1}^{m}\be_i\in H_+(X).$$
With notation as in Section \ref{g0config_subs}, let
$$\cM_{\un\be}=\M_{0,\eset}^*(X,\un\be) \qquad\hbox{and}\qquad
\ov\cM_{\un\be}=\ov\M_{0,\eset}^*(X,\un\be).$$
We denote by $\cM_{\un\be,1}$ and $\ov\cM_{\un\be,1}$ the spaces of 
pairs $(\cC,x)$ such that $\cC\!\in\!\cM_{\un\be}$ and $x\!\in\!\cC$ is 
a nonsingular point of $\cC$ in the first case and 
$\cC\!\in\!\ov\cM_{\un\be}$ and $x\!\in\!\cC$ is any point of $\cC$ 
in the second case.

Let $\cS\!\lra\!\cM_{\un\be}$ be a family of deformations in $X$
of curves in $\cM_{\un\be}$.
In other words, the fiber $\cS_{\cC}$ of $\cS$ over $\cC\!\in\!\cM_{\un\be}$
contains $\cC$ and
$$\dim\, \cS_{\cC}=\dim\,\cM_{|\un\be|,1}-\dim\,\cM_{\un\be}=m.$$
There is a fibration 
\BE{Sfibr_e}
\pi_{\cC}:\cS_{\cC}\lra\De\!\subset\!\C^{m-1}\EE
giving the universal family of deformations of $\cC$.
If $m\!=\!1$, then $\cS\!=\!\cM_{|\un\be|,1}$. 
If $m\!=\!3$, $\cS$ is a small neighborhood
of $\cM_{\un\be,1}$ in~$\ov\cM_{|\un\be|,1}$.

If $\ev\!:\cM_{\un\be,1}\!\lra\!X$ is the evaluation map at the marked point, 
the bundle
\BE{Qdfn_e}Q =\ev^*TX\big/T\cS\lra\cM_{\un\be,1}\EE
extends naturally over $\ov\cM_{\un\be,1}$ so that 
there is an exact sequence
\BE{Qses_e}0\lra  f^*T\ov\cM_{\un\be}\lra Q\lra \N_{|\un\be|}\lra 0,\EE
where $f\!:\ov\cM_{\un\be,1}\!\lra\ov\cM_{\un\be}$ is the forgetful map
and $\N_{|\un\be|}$ is the normal bundle to the family of simple curves
of class~$|\un\be|$.

Similarly to Section 3.3 in \cite{LZ}, we choose a family of ``exponential'' 
maps
\BE{expprp_e}
\exp^{\cC}\!: TX\lra X \qquad\hbox{such that}\qquad
\exp^{\cC}_x(v)\in \cS_{\cC} ~~\hbox{if}~~
x\in\cC,\,v\in T_x\cS_{\cC},\,|v|<\de(\cC),\EE
for some $\de\!\in\!C^{\i}(\cM_{\bar\Ga};\R^+)$.
Below we will place additional assumptions on $\exp^{\cC}$ as needed.\\

For an ideal Calabi-Yau $n$-fold with $n\!\ge\!6$, the above stratification
would need to be refined further based on the deviation of the normal bundles
of curves in $\cM_{\un\be,1}$ from balanced splitting.
The arguments in Sections \ref{g1ghost_subs}-\ref{g1eff_subs}
below apply to the strata with balanced splitting with minor changes.
The main change here is that  the map $\ev$ is no longer an immersion, and
one would need to pass to a blowup of $\ov\cM_{\un\be,1}$ to obtain 
analogues of the vector bundle $Q$ and the short exact sequence~\e_ref{expprp_e}. 
The strata with unbalanced splittings need to be treated separately,
with the conclusion that they do not contribute to the genus $1$ 
Gromov-Witten invariants
under certain assumptions on~$X$.

\subsection{Strata with ghost principal component I}
\label{g1ghost_subs}

Here we describe the contribution to $N_{1,*}$ from a stratum $\cZ_{\T}$ 
consisting of maps  $u\!:\Si\lra\!X$ that are constant on the principal, 
genus-carrying, component(s) $\Si_P$ of $\Si$.
We show  $\cZ_{\T}$ does not contribute to~$N_{1,*}$
unless $\cZ_{\T}$ is of type \ref{basicghost_item}.

For each $m\!\in\!\Z^+$, let
$\ov\cM_{1,m}$ be the moduli space 
of stable curves of genus 1 with $m$ marked points. 
Let $\E\!\lra\!\ov\cM_{1,m}$ be the Hodge line bundle of holomorphic
differentials.
For each $i\!\in\![m]$, denote by $L_i\!\lra\!\ov\cM_{1,m}$
the universal tangent line bundle at the $i^{th}$ marked point.
Let
$$s_i\in\Ga\big(\ov\cM_{1,m},\Hom(L_i,\E^*)\big)$$
be the homomorphism induced by the natural pairing of 
tangent and cotangent vectors at the $i$th marked point.
Denote by 
$$\cM_{1,m},\cM_{1,m}^{\eff}\subset\ov\cM_{1,m}$$
the subspaces consisting of nonsingular
curves and of curves $\cC$ with 
no bubble components ($\cC$ is either a nonsingular
 genus 1 curve or is a circle  of rational curves).

Let $L_1\!\lra\!\ov\M_{0,1}(X,\be)$
be the universal tangent line bundle at the marked point.
Denote by
$$\cD_1\in\Ga\big(\ov\M_{0,1}(X,\be),\Hom(L_1,\ev_1^*TX)\big)$$
the natural homomorphism induced by the derivative of the map
at the marked point.
For $m\!\in\!\Z^+$, let 
\begin{equation*}\begin{split}
\ov\M_{(0,m)}(X,\be)=\big\{
(b_i)_{i\in[m]}\!\in\!\prod_{i=1}^{m}\ov\M_{0,\{0\}}(X,\be_i)\!:\, &
\be_i\!\in\!H_+(X),\, \sum_{i=1}^{m}\be_i\!=\!\be,\\
&\ev_0(b_i)\!=\!\ev_0(b_{i'})\, \forall i,i'\!\in\![m]\big\}.
\end{split}\end{equation*}
There is a well-defined evaluation map
$$\ev_0\!: \ov\M_{(0,m)}(X,\be)\lra X, \qquad 
(b_i)_{i\in[m]}\lra\ev_0(b_i),$$
which is independent of the choice of $i$.
Let
$$\pi_i\!: \ov\M_{(0,m)}(X,\be)\lra 
\bigsqcup_{\be_i\in H_+(X)}\ov\M_{0,\{0\}}(X,\be_i)$$
be the projection onto the $i^{th}$ component.
Denote by
$$\M_{(0,m)}^{\eff}(X,\be)\subset\ov\M_{(0,m)}(X,\be)$$
the subset consisting of the tuples  $(u_i)_{i\in[m]}$  
such that for each $i\!\in\![m]$ the restriction of $u_i$ to
the domain component carrying the marked point $0$ 
is not constant.

The stratum $\cZ_{\T}$ admits a decomposition
\BE{Zdecomp_e} \cZ_{\T}=
\big(\cZ_{\T,P}\times\cZ_{\T,PB}\times\cZ_{\T,B}\big)\big/S_{m_B}, \EE
where $\cZ_{\T,P}$ is a stratum of $\cM_{1,m_P}^{\eff}$ for 
some $m_P\!\in\!\Z^+$,
$\cZ_{\T,B}$ is a stratum of $\M_{(0,m_B)}^{\eff}(X,\be)$ for some $m_B\!\in\!\Z^+$,
and $\cZ_{\T,PB}$ is a product of moduli spaces of 
irreducible stable genus 0 curves.
The stratum $\cZ_{\T,P}$ consists of curves of a fixed topological type,
while the elements of $\cZ_{\T,B}$ are tuples of stable maps from domains 
of fixed topological types so that the image of the restriction 
of the map to each component is of a specified homology class
and multiplicity.
The requirement that  
$$\cZ_{\T,P}\subset\cM_{1,m_P}^{\eff} \qquad\hbox{and}\qquad
\cZ_{\T,B}\subset\M_{(0,m_B)}^{\eff}(X,\be)$$ 
implies that the decomposition \e_ref{Zdecomp_e} is well-defined.
Let
$$\pi_P,\pi_B\!: \cZ_{\T,P}\times\cZ_{\T,PB}\times\cZ_{\T,B}
\lra \cZ_{\T,P},\cZ_{\T,B}$$
denote the projection maps. The quotient is by the 
automorphism groups $S_{m_B}$ of the data.

If $X$ is an ideal CY $5$-fold, $\cZ_{\T,B}$ is smooth.
The cokernels of the linearizations $D_b$ of the $\dbar$-operator 
along $\cZ_{\T}$ form the obstruction bundle
\BE{Odecomp_e} \fO=\fO_{PB}\oplus\pi_B^*\fO_B =
\pi_P^*\E^*\!\otimes\!\pi_B^*\ev_0^*TX\oplus\pi_B^*\fO_B,\EE
where $\fO_B\!\lra\!\cZ_{\T,B}$ is the obstruction bundle
associated with the moduli space $\ov\M_{(0,m)}(X,\be)$.
Let $\bar\nu\!\in\!\Ga(\cZ_{\T},\fO)$ be the section induced by $\nu$:
$\bar\nu(b)$ is the projection of $\nu(b)$ to the cokernel of~$D_b$.
We write 
$$\bar\nu_{PB},\bar\nu_B\in \Ga\big(\cZ_{\T},\fO_{PB}\big),
\Ga\big(\cZ_{\T},\fO_B\big)$$
for the two components of~$\bar\nu$.

There is a natural projection map
$$\bar\pi:\cZ_{\T,B}\lra\cM_{\un\be,1},$$
sending $\pi_B([\Si,u])$ to $(u(\Si),u(\Si_P))$.
Denote by 
$$\bar\nu_{PB}^{\perp}\in
\Ga\big(\cZ_{\T},\pi_P^*\E^*\!\otimes\!\pi_B^*\bar\pi^*Q\big)$$
the image of $\bar\nu_{PB}$ under the natural projection map.
Let 
$$f\!:\ov\cM_{1,m_P}\lra\ov\cM_{1,1}$$
be the forgetful map, dropping all but the first marked point.
The restriction of the bundle
\BE{partcok_e}
\pi_1^*\E^*\otimes \pi_2^*Q\lra\ov\cM_{1,1}\times\ov\cM_{\un\be,1}\EE
to any boundary stratum $\cZ_{\Ga}$ contains a subbundle $\fO_{\Ga}$
such that 
\BE{rkmindim_e}
\rk\,\fO_{\Ga}-\dim\,\cZ_{\Ga}>
\rk\,\big(\pi_1^*\E^*\!\otimes\!\pi_2^*Q\big)
-\dim\,\big(\ov\cM_{1,1}\!\times\!\ov\cM_{\un\be,1}\big)=0
\EE
and $\big\{(f\!\circ\!\pi_P)\times(\bar\pi\!\circ\!\pi_B)\big\}^*\fO_{\Ga}$
is a quotient of the cokernel bundle over a boundary stratum of~$\ov\cZ_{\T}$.
Thus, we can choose a section $\bar\nu_{\un\be}$ of \e_ref{partcok_e}
with  all zeros  transverse and contained in 
$\cM_{1,1}\!\times\!\cM_{\un\be,1}$ and such that there exists $\nu$ as above satisfying
$$\bar\nu_{PB}^{\perp}=
\big\{(f\!\circ\!\pi_P)\times(\bar\pi\!\circ\!\pi_B)\big\}^*
\bar\nu_{\un\be}.$$

It is shown in the next section that  the contribution of $\cZ_{\T}$ to 
$N_{1,*}$ comes from $\bar\nu_{PB}^{\perp\,-1}(0)$.
Thus, if $\cZ_{\T,P}\!\not\subset\!\cM_{1,m_P}$, then 
$\bar\nu_{PB}^{\perp\,-1}(0)$ is empty for a good choice of $\nu$ 
by \e_ref{rkmindim_e} and 
the stratum $\cZ_{\T}$ does not contribute to $N_{1,*}$.
Otherwise, $\bar\nu_{PB}^{\perp\,-1}(0)$ is the preimage of a finite 
subset in~$\cM_{1,1}\!\times\!\cM_{\un\be,1}$.
It decomposes into connected components
\BE{Zdecomp_e5}\bar\nu_{PB}^{\perp\,-1}(0)
=\bigsqcup_{(\cC,x)\in\pi_2(\bar\nu_{\un\be}^{-1}(0))}
\!\!\!\!\!\!\!\cZ_{\cC,x},\EE
where $\cC$ is a $\un\be$-curve and $x$ is a nonsingular point of $\cC$.
Then, $\tC_{\T}(\be)$ is the number of zeros of a map $\vph_{t\nu}$
from  the vector bundle $F$ of gluing parameters to~$\fO$ over
each of the components~$\cZ_{\cC,x}$. 
The projection of $\vph_{t\nu}$ in the decomposition~\e_ref{Odecomp_e} 
onto $\pi_P^*\E\!\otimes\!T_x\cS_{\cC}/T_x\cC$ is essentially the same 
as the projection of $t\nu$, which we denote by~$t\ti\nu$.
Since $\ti\nu$ is a section of a trivial bundle over~$\cZ_{\cC,x}$,
it can be chosen not to vanish if $m\!>\!1$.
Thus, $\tC_{\T}(\be)\!=\!0$ if $m\!>\!1$.
On the other hand, the second component of $\vph_{t\nu}$ with respect to 
the decomposition~\e_ref{Odecomp_e} is essentially~$t\bar\nu_B$.
It does not vanish on $\bar\nu_{PB}^{\perp\,-1}(0)$ for dimensional reasons
if $m\!=\!1$, but $|\Ver|\!>\!2$.
Thus,  $\tC_{\T}(\be)\!=\!0$ if $\T$ is not basic.

Finally, if $\T$ is basic, the principal component of every element of 
$\cZ_{\cC,x}$ is a fixed nonsingular genus 1 curve $\Si_P$
with one special point $z_1$ and  
$$\cZ_{\cC,x}\approx \M_{0,1}(\P_p,d),$$
where $\M_{0,1}(\P_p,d)\subset\ov\M_{0,1}(\P,d)$
is the subspace of elements $[\Si,u]$ such that 
$\Si$ is nonsingular and $\ev_1([\Si,u])\!=\!p$ for a fixed $p\!\in\!\P$.
Let
$$\D\in\Ga\big(\ov\cM_{1,1}\!\times\!\ov\M_{0,1}(\P,d),
\Hom(\pi_1^*L_1\!\otimes\!\pi_2^*L_1,\pi_1^*\E^*\!\otimes\!\pi_2^*\ev_1^*T\P)\big)$$
be given by
\BE{Ddnf_e} \D(v\!\otimes\!w)=s_1(v)\otimes \cD_1(w) \,.\EE
The first component of $\vph_{t\nu}$ with respect to 
the decomposition~\e_ref{Odecomp_e} is essentially 
\BE{Omap_e1}
F=\pi_1^*L_1|_{z_1}\!\otimes\!\pi_2^*L_1\lra 
\fO_{PB}=\E_{\Si_P}^*\!\otimes\!T_x\cC,\qquad
\ups\lra \D(\ups)+t\bar\nu_{PB}.\EE
Let $\U$ be the universal curve over $\ov\M_{0,1}(\P,d)$,
with structure map~$\pi$ and evaluation map~$\ev$:
\BE{univcurv_e}\xymatrix{\U \ar[d]^{\pi} \ar[r]^{\ev} & \P \\
\ov\M_{0,1}(\P,d).}\EE
The restriction of $\bar\nu_B$ to $\bar\nu_{PB}^{\perp\,-1}(0)$ is 
a section of 
$$\fO_B=R^1\pi_*\ev_*\big(\cO(-1)\oplus\cO(-1)\big)
\lra \ov\M_{0,1}(\P_p,d).$$
Thus, by the Aspinwall-Morrison and divisor formulas, 
as in Section 1.1 in~\cite{KlP},
\BE{AM_e1} ^{\pm}\big|\bar\nu_B^{-1}(0)\big| = \frac{1}{d^2}.\EE
On the other hand, $\bar\nu_{PB}^{\perp}$ is a section of
$$\pi_1^*\E^*\otimes\pi_2^*\big(f^*T\ov\cM_{\be}\oplus\N_{\be}\big)
\lra \ov\cM_{1,1}\times\ov\cM_{\be,1},$$
see \e_ref{Qdfn_e}.
Therefore,
\BE{contr_e1} 
^{\pm}\big|\bar\nu_{PB}^{\perp\,-1}(0)\big|=
-\frac{1}{24}\int_{\ov\cM_{\be,1}}\big(
f^*c_2(\ov\cM_{\be})\ev_1^*c_1(\N_{\be})
+f^*c_1(\ov\cM_{\be})\,c_2(\N_{\be})\big).\EE
Since \e_ref{Omap_e1} has a unique zero in every fiber of $F$ over
$\bar\nu_{PB}^{\perp\,-1}(0)\!\cap\!\bar\nu_B^{-1}(0)$ and the restriction
of $\N_{\be}$ to a fiber of $f$ is of degree $-2$, 
equations \e_ref{AM_e1} and \e_ref{contr_e1} imply 
\BE{ghcontr_e}
\tC_{\T_{\gh}(\be,d)}(d\be)=
\frac{1}{24d^2}\bigg(\int_{\ov\cM_{\be}}2c_2(\ov\cM_{\be})
-\int_{\ov\cM_{\be,1}} f^*c_1(\ov\cM_{\be})c_2(\N_{\be})\bigg).\EE
We have proved the second scaling identity in \e_ref{g1thm_e}.
The equation
\BE{ghx}
\tC_{\T_{\gh}(\be,1)}(\be)=
\frac{1}{24} \int_{\ov\cM_{\be}}2c_2(\ov\cM_{\be})-c_1^2(\ov\cM_{\be})\EE
is obtained from \eqref{ghcontr_e} from relations
\eqref{cN_e} and \eqref{cF_e}.

\subsection{Strata with ghost principal component II}
\label{g1ghost_subs2}

\noindent
We continue with the setup of Section \ref{g1ghost_subs}.
For each $[\Si_u,u]\!\in\!\cZ_{\T}$, denote by $\Si_u^0\!\subset\!\Si_u$ 
the largest union of irreducible components of $\Si_u$ that 
contains the principal component(s) of~$\Si_u$ and on which~$u$ is constant.
The topological types of $\Si_u$ and~$\Si_u^0$ are independent of
the choice of $[\Si_u,u]\!\in\!\cZ_{\T}$.

The bundle of gluing parameters (or smoothing of the nodes) over 
$F\!\lra\!\cZ_{\T}$ is a direct sum of line bundles (up to a quotient 
by a finite group).
Let $F^{\eset}\!\subset\!F$ be the subspace of smoothings with all 
components nonzero, smoothings that do not leave any nodes.
If $\ups\!\in\!F_u^{\eset}$ is sufficiently small, there is a $C^\infty$-map
$$q_{\ups}\!: \Si_{\ups}\lra\Si_u,$$
where $\Si_u$ is the domain of $u$ and $\Si_{\ups}$ is a genus 1 
Riemann surface with thin necks replacing the nodes of $\Si_u$,
see Section 2.2 of \cite{gluing}.
This map determines Riemannian metrics and weights on $\Si_{\ups}$ 
which induce the $L^p_1$- and $L^p$ Sobolev norms, 
$\|\cdots\|_{\ups,p,1}$ and  $\|\cdots\|_{\ups,p}$, with $p\!>\!2$,
appearing below,  see Section 3.3 in \cite{gluing}.
These norms are equivalent to the ones used in Section 3 of \cite{LT}.
Let $\Si_{\ups}^0\!=\!q_{\ups}^{-1}(\Si_u^0)$.

We take the approximately holomorphic map corresponding to $\ups\!\in\!F_u$
 to be
$$u_{\ups}\!=\!u\!\circ\!q_{\ups}\!: \Si_{\ups}\lra X.$$ 
The map satisfies
\BE{preglmap_e} \big\|\dbar u_{\ups}\big\|_{\ups,p}\le C(u)|\ups|^{1/p}.\EE
Let 
$$D_{\ups}\!: \Ga(\ups)\!=\!\Ga(\Si_{\ups},u_{\ups}^*TX) 
\lra \Ga^{0,1}(\ups)\!=\!
\Ga(\Si_{\ups},T^*\Si_{\ups}^{0,1}\!\otimes\!u_{\ups}^*TX)$$
be the linearization of the $\dbar$-operator at $u_{\ups}$
defined using the Levi-Civita connection of a Kahler metric $g_{X,u}$ on $X$.
As in Sections~2 and~4.1 in \cite{g1comp}, 
we can construct splittings
\BE{Gadecomp_e} \Ga(\ups)=\Ga_-(\ups)\oplus\Ga_+(\ups) \qquad\hbox{and}\qquad
\Ga^{0,1}(\ups)=\Ga^{0,1}_{-;PB}(\ups)\oplus\Ga^{0,1}_{-;B}(\ups)
\oplus\Ga^{0,1}_+(\ups),\EE
and isomorphisms 
\BE{Rmap_e} R_{\ups}\!: \fO_{PB}\oplus\pi_B^*\fO_B\lra 
\Ga^{0,1}_{-;PB}(\ups)\oplus\Ga^{0,1}_{-;B}(\ups)\EE
with the following properties:
\begin{enumerate}[label=(G\arabic*)]
\item\label{Disom_item} $D_{\ups}\!:\Ga_+(\ups)\lra\Ga^{0,1}_+(\ups)$ 
is an isomorphism with the norm of the inverse bounded 
independently of $\ups\!\in\!F_u^{\eset}$ (but depending on~$[\Si,u]$),
\item\label{supp_item} the elements of $\Ga^{0,1}_{-;PB}(\ups)$ are supported 
on a small neighborhood of~$\Si_{\ups}^0$,
\item\label{proj_item} 
if $\pi^{0,1}_{-;PB}\!:\Ga^{0,1}(\ups)\lra\Ga^{0,1}_{-;PB}(\ups)$
is the projection in the second decomposition~\e_ref{Gadecomp_e},
\BE{minuspr_e} 
\big\|\pi^{0,1}_{-;PB}D_{\ups}\xi\big\|_{\ups,2}
\le C(u)|\ups|\|\xi\|_{\ups,p,1}
 ~~~\forall\,\xi\in\Ga(\ups),\EE
\item\label{proj_item2}  if $|\Ver|\!=\!2$ (and thus $\cZ_{\T}$ is basic),
\BE{minuspr_e2} 
\pi^{0,1}_{-;PB}\dbar u_{\ups}=R_{\ups}\D\ups,\EE
with $\D$ as in \e_ref{Ddnf_e},
\item\label{neighb_item} every map $\ti{u}\!:\Si\!\lra\!X$, 
where $\Si$ is a smooth genus-one Riemann surface, 
that lies in a small neighborhood of~$\cZ_{\T}$ can be written
uniquely~as $\ti{u}\!=\!\exp_{u_{\ups}}\xi$ for small $\ups\!\in\!F^{\eset}$
and $\xi\!\in\!\Ga_+(\ups)$.
\end{enumerate}
Let 
$$\pi^{0,1}_+\!:\Ga^{0,1}(\ups)\lra\Ga^{0,1}_+(\ups)
\qquad\hbox{and}\qquad
\pi^{0,1}_{-;B}\!:\Ga^{0,1}(\ups)\lra\Ga^{0,1}_{-;B}(\ups)$$
be the component projections  in the second decomposition~\e_ref{Gadecomp_e}.

The relation \e_ref{CR_e2} for $\ti{u}\!=\!\exp_{u_{\ups}}\xi$ is equivalent~to
\BE{CR_e3} \dbar u_{\ups}+D_{\ups}\xi+t\nu_{\ups}
+N_{\ups}(\xi)+tN_{\nu,\ups}(\xi)=0,\EE
with $N_{\ups}$ and $N_{\nu,\ups}$ satisfying
\begin{equation}\label{Nterm_e}\begin{split}
\big\|N_{\ups}(\xi)-N_{\ups}(\xi')\big\|_{\ups,p} &\le 
C(u)\big(\|\xi\|_{\ups,p,1}\!+\!\|\xi'\|_{\ups,p,1}\big)
\big\|\xi\!-\!\xi'\big\|_{\ups,p,1},\\
\big\|N_{\nu,\ups}(\xi)-N_{\nu,\ups}(\xi')\big\|_{\ups,p} &\le 
C(u)\big\|\xi\!-\!\xi'\big\|_{\ups,p,1},
\end{split}\end{equation}
if $\ups\!\in\!F_u^{\eset}$.
For a good choice of identifications,
\BE{Nterm_e0} \pi_{-;PB}^{0,1}N_{\ups}\xi=0 \qquad~\forall\xi\in\Ga(\ups).\EE
By the Contraction Principle and \ref{Disom_item}, the equation
$$\pi_+^{0,1}\big(\dbar u_{\ups}+D_{\ups}\xi+t\nu_{\ups}
+N_{\ups}(\xi)+tN_{\nu,\ups}(\xi)\big)=0$$
has a unique small solution $\xi_{t\nu}(\ups)\!\in\!\Ga_+(\ups)$.
By \e_ref{preglmap_e}, it satisfies
\BE{xinorm_e} \big\|\xi_{t\nu}(\ups)\big\|_{\ups,p,1}
\le C(u)\big(|\ups|^{1/p}\!+\!t\big).\EE
Thus, the number of solutions of \e_ref{CR_e2} near $\cZ_{\T}$
is the number of solutions of the equation
\BE{CR_e4} \dbar u_{\ups}+D_{\ups}\xi_{t\nu}(\ups)+t\nu_{\ups}
+N_{\ups}\big(\xi_{t\nu}(\ups)\big)+
tN_{\nu,\ups}\big(\xi_{t\nu}(\ups)\big)=0
\in \Ga^{0,1}_{-;PB}(\ups)\oplus\Ga^{0,1}_{-;B}(\ups).\EE
This is an equation on $\ups\!\in\!F_u^{\eset}$ with $|\ups|\!<\!\de(u)$
for some $\de\!\in\!C^{\i}(\cZ_{\T};\R^+)$.

For each $[\Si,u]\!\in\!\cZ_{\T}$, 
let $\cC_u\!=\!u(\Si)$ and $\cS_u\!=\!\cS_{\cC_u}$,
see the end of Section \ref{g1summ_subs}.
The pregluing map $u_{\ups}$ satisfies
$$u_{\ups}(\Si_{\ups})\subset\cC_u\subset\cS_u.$$
We can choose the splittings \e_ref{Gadecomp_e} so that 
they restrict to splittings for vector fields and $(0,1)$-forms
along $u_{\ups}$ with values in $T\cS_u$ and
\ref{Disom_item} holds when restricted to $T\cS_u$.
If $\exp_{u_{\ups}}\xi$ is defined using the ``exponential''  $\exp^{\cC_u}$,
the operators $D_{\ups}$ and $N_{\ups}$ in~\e_ref{CR_e3} preserve 
$T\cS_u$ as well.
Therefore, 
\BE{xi0term_e} \xi_0(\ups)\in\Ga(\Si_{\ups},u_{\ups}^*T\cS_u), \qquad
D_{\ups}\xi_0(\ups),N_{\ups}\big(\xi_0(\ups)\big)
\in \Ga\big(\Si_{\ups},T^*\Si_{\ups}^{0,1}\!\otimes\!u_{\ups}^*T\cS_u\big).\EE
On the other hand, by~\ref{Disom_item} and~\e_ref{Nterm_e},
\BE{xidiff_e} 
\big\|\xi_{t\nu}(\ups)-\xi_0(\ups)\big\|_{\ups,p,1}\le C(u)t,\EE
if $\ups\!\in\!F_u$ is sufficiently small.
Taking the projection $\pi_{-;PB}^{\perp}$ of~\e_ref{CR_e4}
to $\pi_P^*\!\E^*\!\otimes\!\bar\pi^*Q$,
we thus find that any solution $\ups\!\in\!F_u$ of~\e_ref{CR_e4}
satisfies
$$\big\|\bar\nu_{PB}^{\perp}(u)\big\|\le \ve(t,\ups)$$
for some function $\ve\!:\R\!\times\!F^{\eset}\!\lra\!\R^+$ 
approaching~$0$ as $(t,\ups)$ approaches~$0$.
Therefore, all solutions of~\e_ref{CR_e2} lie in a small neighborhood
of $\bar\nu_{PB}^{\perp\,-1}(0)\!\subset\!\cZ_{\T}$, 
as claimed in Section \ref{g1ghost_subs}.

If $m\!=\!1$, for a good choice of $R_{\ups}$ on $\pi_B^*\fO_B$
\BE{CR_e5} 
\bllrr{\eta,\eta_-}_{\ups,2}=0 \qquad\forall~
\eta\!\in\!\Ga\big(\Si_{\ups},
T^*\Si_{\ups}^{0,1}\!\otimes\!u_{\ups}^*T\cS_u\big),
\, \eta_-\!\in\!\Ga^{0,1}_{-;B}(\ups). \,\footnotemark\EE
Taking the projection of \e_ref{CR_e4} onto $\Ga^{0,1}_{-;B}(\ups)$
and using~\e_ref{Nterm_e}, \e_ref{xi0term_e}, \e_ref{xidiff_e}, and~\e_ref{CR_e5},
\footnotetext{In this case, $\fO_B$ is isomorphic to the cokernel 
of a $\dbar$-operator on $\N_{|\un\be|}$.}
we obtain
\BE{CR_e6} 
t\bar\nu_B(u)+t\eta(t,\ups)=0,\EE
for some $\eta(t,\ups)$ approaching~$0$ as $(t,\ups)$ approaches~$0$.
Since $u$ is $d\!:\!1$ cover of the smooth curve $\cC_u$,
the dimension of the projection of $\bar\nu_{PB}^{\perp\,-1}(0)$
onto the third component in the decomposition \e_ref{Zdecomp_e} is of dimension
at most $2d\!-\!2m_B$.\footnote{This is the 
dimension of the space of degree-$d$ 
covers of $\bP^1$ by $m_B$ copies of $\bP^1$.
The dimension is less than $2d\!-\!2m_B$ unless  
$\cZ_{\T,B}$ is the main stratum of $\ov\M_{(0,m_B)}(X,\be)$.}
Since the rank of $\fO_B$ is $2d\!-\!2$, \e_ref{CR_e6} has no solutions 
for a generic choice of~$\nu$ unless $\cZ_{\T}$ is  described 
by \ref{basicghost_item} of Section \ref{g1dfn_subs}.
If $\cZ_{\T}$ is of type~\ref{basicghost_item}, the number of solutions
of~\e_ref{CR_e4} is the same as the number of small solutions of
\begin{gather}\label{CR_e7}
\D\ups+t\bar\nu_{PB}(u)+\eta(t,\ups)=0\in
\E^*\!\otimes\!T_{u(\Si_P)}\cC_u\\
\ups\in \pi_1^*L_1\!\otimes\!\pi_2^*L_1|_u,\quad
u\in\bar\nu_{PB}^{\perp\,-1}(0)\!\cap\!\bar\nu_B^{-1}(0)
\subset\cM_{1,1}\!\times\!\M_{0,1}(\bP^1;d),\notag
\end{gather}
with the error term $\eta(t,\ups)$ satisfying
$$\big\|\eta(t,\ups)\big\|_{\ups,2}\le 
\ve(t,\ups)\big(t\!+\!|\ups|);$$
see \e_ref{minuspr_e}, \e_ref{minuspr_e2}, \e_ref{Nterm_e}, and~\e_ref{xinorm_e}. 
If $\nu$ is generic, $\cD_1$ and thus $\D$ are nowhere zero on 
the finite set $\bar\nu_{PB}^{\perp\,-1}(0)\!\cap\!\bar\nu_B^{-1}(0)$.
By the same rescaling and cobordism argument as in Section 3.1 
of~\cite{g2n2and3},
the number of small solutions of~\e_ref{CR_e7} is the same as the number of 
solutions of 
$$\D\ups+\bar\nu_{PB}(u)=0, \qquad
\ups\in \pi_1^*L_1\!\otimes\!\pi_2^*L_1
\big|_{\bar\nu_{PB}^{\perp\,-1}(0)\cap\bar\nu_B^{-1}(0)}.$$
There is one solution for each of the elements of 
$\bar\nu_{PB}^{\perp\,-1}(0)\cap\bar\nu_B^{-1}(0)$.
This concludes the consideration of the $m\!=\!1$ case.

We will next show that \e_ref{CR_e4} has no solution if $m\!>\!1$. 
Let $\cZ_{\cC,x}$ be as in~\e_ref{Zdecomp_e5}.
Since $x$ is a nonsingular point of $\cC$,
on a neighborhood~$U$ of~$x$ in~$\cC$ there is an orthogonal decomposition
\BE{TSdecomp_e}TS_{\cC}|_{U}=T\cC|_{U}\oplus\N\cC|_{U}.\EE
We can assume that the ``exponential" map $\exp^{\cC}$ satisfies
\BE{expCprp_e2}
\pi_{\cC}\big(\exp_y^{\cC}v\big)=d\pi_{\cC}|_xv
\qquad\forall\,y\in U,\,v\in T_y\cS_{\cC},\,|v|<\de, \EE
with $\pi_{\cC}$ as in~\e_ref{Sfibr_e}.
For any $[\Si,u]\!\in\!\cZ_{\T}$ in a small neighborhood of $\cZ_{\cC,x}$,
let $W_u\!=\!u^{-1}(U)$ be an open neighborhood of $\Si_u^0$ in $\Si_0$.
We can assume that 
every element $\eta$ of $\Ga_{-;PB}^{0,1}(\ups)$
is supported in $W_{\ups}\!=\!q_{\ups}^{-1}(W_u)$,
whenever $\ups\!\in\!F_u^{\eset}$ is sufficiently small.
With $D_{\ups}^*$ denoting the formal adjoint of $D_{\ups}$ with respect
to the inner-product $\llrr{\cdot,\cdot}_{\ups,2}$, let 
\begin{equation}\label{CR_9}\begin{split}
\Ga_{+-}(\ups)&=\big\{\xi\!\in\!\Ga_+(\ups)\!: 
\bllrr{\xi,D_{\ups}^*R_{\ups}\eta}\!=\!0~\forall\,
\eta\!\in\!\E_{\Si_P}^*\!\otimes\!\N_{u(\Si_P)}\cC_u\big\}
\subset\Ga(\ups),\\
\Ga_{++}(\ups)&=\big\{D_{\ups}^*R_{\ups}\eta\!:
\eta\!\in\!\E_{\Si_P}^*\!\otimes\!\N_{u(\Si_P)}\cC_u\big\}
\subset D_{\ups}^*\Ga^{0,1}_{-;PB}(\ups)\subset\Ga(\ups).
\end{split}\end{equation}
An explicit expression for $D_{\ups}^*R_{\ups}\eta$ is given 
in the proof of Lemma 2.2 in \cite{g2n2and3}.
Section 2.3 of \cite{g2n2and3} implies that we can take
\BE{Gadecomp_e2}
\Ga_+(\ups)=\Ga_{++}(\ups)\oplus \Ga_{+-}(\ups).\EE
In particular, the proof of Lemma 2.6 shows that the limits of 
the spaces $\Ga_{++}(\ups)$ as $\ups\!\lra\!0$ are 
orthogonal to the limits of the spaces~$\Ga_-(\ups)$.
The decomposition~\e_ref{Gadecomp_e2} is $L^2$-orthogonal by~\e_ref{CR_9}
and
\BE{Gaplusplus_e}\|\xi\|_{\ups,p,1}\le C(u)\|\xi\|_{\ups,2}
\qquad\forall\xi\in\Ga_{++}(\ups),\EE
see the proof of Lemma~2.2 in~\cite{g2n2and3}.

Let $\xi_{t\nu}^+(\ups)$ and $\xi_{t\nu}^-(\ups)$ be the components of 
$\xi_{t\nu}(\ups)$ with respect to the decomposition~\e_ref{Gadecomp_e2}.
Denote by $\ti\nu(u)\!\in\!\E_{\Si_P}^*\!\otimes\!\N_{u(\Si_P)}\cC_u$ 
the projection of $\nu(u)$ to $\E_{\Si_P}^*\!\otimes\!\N_{u(\Si_P)}\cC_u$.
Since $\ti\nu$ is a section of a trivial bundle near~$\cZ_{\cC,x}$,
we can assume that it has no zeros on $\cZ_{\cC,x}$.
In the next paragraph we will show  
\BE{CR_e8} \big\|\xi_{t\nu}^+(\ups)\big\|_{\ups,p,1}\le C(u)t.\EE
Assuming this is the case, we project both sides of~\e_ref{CR_e4} onto 
$$R_{\ups}\big(\E_{\Si_P}^*\!\otimes\!\N_{u(\Si_P)}\cC_u\big)
\subset\Ga_{-;PB}^{0,1}(\ups)$$
and take the  preimage under $R_{\ups}$.
Since the projections of $\dbar u_{\ups}$ and $N_{\ups}(\xi_{t\nu(\ups)})$
vanish, using the first equation in~\e_ref{CR_9}, \e_ref{CR_e8},
\e_ref{minuspr_e2}, and~\e_ref{Nterm_e}, we obtain
$$t\ti\nu(u)+\eta(t,\ups)=0$$
with $\eta(t,\ups)$ satisfying
$$\big\|\eta(t,\ups)\big\|_{\ups,2}\le \ve(t,u)t.$$
However, this is impossible if $t$ and $\ups$ are sufficiently small
(``small" depending continuously on~$u$), since $\ti\nu$ has no zeros 
over~$\cZ_{\cC,x}$.

We now verify \e_ref{CR_e8}.
Let
$$\ti\xi_{t\nu}(\ups)\!=\!\pi_{\cC_u}\!\circ\!\exp_{u_{\ups}}\xi_{t\nu}(\ups)
\!:\Si_{\ups}\lra\C^{m-1}, \qquad
\ti\xi_{t\nu}^{\pm}(\ups)=d\pi_{\cC_u}\circ\xi_{t\nu}^{\pm}(\ups).$$
Since $\xi_{t\nu}^+(\ups)$ is supported on $W_u$, by~\e_ref{expCprp_e2}
\BE{CR9c}
\blr{\ti\xi_{t\nu}^+(\ups),\ti\xi_{t\nu}^{\pm}(\ups)}_z
=\blr{\xi_{t\nu}^+(\ups),\xi_{t\nu}^{\pm}(\ups)}_z
\qquad\forall\,z\!\in\!\Si_{\ups}\,.\EE
By~\e_ref{expCprp_e2}, we also have
\BE{CR9d}\ti\xi_{\ups}|_{W_{\ups}}=
\ti\xi_{\ups}^+|_{W_{\ups}}+\ti\xi_{\ups}^-|_{W_{\ups}}\,.\EE
Since~\e_ref{CR_e4} is equivalent to~\e_ref{CR_e2} for 
$\ti{u}\!=\!\exp_{u_{\ups}}\xi_{t\nu}(\ups)$,
\BE{CR_e9a} 
\big\|\dbar\ti\xi_{t\nu}(\ups)\big\|_{\ups,p}\le C(u)t.\EE
Since the operator
$$L^p_1\big(\Si_{\ups},\C^{m-1}\big)\lra
L^p\big(\Si_{\ups},T^*\Si^{0,1}_{\ups}\C^{m-1}\big)\oplus\C^{m-1}, \qquad
\ti\xi\lra\bigg(\dbar\ti\xi,\int_{\Si_{\ups}}\ti\xi\,dvol_{\Si_{\ups}}\bigg),$$
is an isomorphism with the norm of the inverse bounded independently of $\ups$
(but depending on~$u$), \e_ref{CR_e9a} implies that 
\BE{CR_e9b}
\big\|\ti\xi_{t\nu}(\ups)-A_{t\nu}(\ups)\big\|_{\ups,p,1}\le C(u)t\EE
for some $A_{t\nu}(\ups)\!\in\!\C^{m-1}$.
Since $\xi_{t\nu}^+(\ups)$ is supported on~$W_u$, 
by~\e_ref{expCprp_e2} and~\e_ref{CR_9},
$$\bllrr{\ti\xi_{t\nu}^+(\ups),A_{t\nu}(\ups)}_{\ups,2}=0.$$
Thus, by~\e_ref{CR9c}, \e_ref{CR9d}, and~\e_ref{CR_e9b},
\begin{equation*}\begin{split}
\big\|\xi_{t\nu}^+(\ups)\big\|_{\ups,2}
=\big\|\ti\xi_{t\nu}^+(\ups)|_{W_{\ups}}\big\|_{\ups,2}
&\le \big\|(\ti\xi_{t\nu}(\ups)\!-\!A_{t\nu}(\ups))|_{W_{\ups}}\big\|_{\ups,2}\\
&\le \big\|\ti\xi_{t\nu}(\ups)-A_{t\nu}(\ups)\big\|_{\ups,p,1}\le C'(u)t.
\end{split}\end{equation*}
The estimate~\e_ref{CR_e8} now follows from~\e_ref{Gaplusplus_e}.

Finally, we comment on the choices made in \e_ref{Gadecomp_e}
and~\e_ref{Rmap_e}.
Choosing the splittings~\e_ref{Gadecomp_e} so that 
\ref{Disom_item} and~\ref{neighb_item} hold is essentially equivalent 
to choosing approximate kernel and cokernel for $D_{\ups}$
that vary smoothly with~$\ups$.
This is easily accomplished in many possible ways, 
including via the construction in Section~3 of~\cite{LT}.
In order to ensure that \ref{supp_item}-\ref{proj_item2} hold, 
$R_{\ups}$ on $\fO_{PB}$ is constructed by pushing harmonic forms on~$\Si_P$ 
over a small neighborhood of~$\Si_{\ups}^0$,
see Section 2.2 of~\cite{g2n2and3}.
Finally, in order to obtain \e_ref{Nterm_e0}, define $\exp_{u_{\ups}}$
and parallel transport using a Kahler metric which is flat near
$u(\Si_P)$, as in Section 2.1 of~\cite{g2n2and3}.

\subsection{Strata with effective principal component}
\label{g1eff_subs}

We determine here the contribution to $N_{1,*}$ from 
a stratum $\cZ_{\T}$ consisting of maps 
$u\!:\Si\lra\!X$ that are not constant on the principal, genus-carrying,
component(s) $\Si_P$ of~$\Si$.
We show that $\cZ_{\T}$ does not contribute to~$N_{1,*}$
unless $\cZ_{\T}$ is of type~\ref{basiceff_item}.

Let $\fO\!\lra\!\cZ_{\T}$ and $F\!\lra\!\cZ_{\T}$ be the obstruction bundle 
and the bundle of gluing parameters as before.
The projection map $\ev^*TX\!\lra\!Q$ induces a surjective homomorphism
\BE{fOsurj_e}\pi^{\perp}\!:\fO\lra\fO^{\perp},\EE
where $\fO^{\perp}|_{[\Si,u]}$ is the cokernel of 
the $\dbar$-operator on $Q$ induced by the $\dbar$-operator~$D_b$ on~$TX$.
By a gluing and obstruction bundle analysis similar to Section~\ref{g1ghost_subs2},
$\tC_{\T}(*)$ is the number of zeros of a bundle map
$$\vph_{t\nu}\!:F\lra\fO$$ 
over $\cZ_{\T}$ for $t$ sufficiently small.
As in the previous case, all zeros of $\vph_{t\nu}$ arise from the zeros~of 
$$\bar\nu^{\perp}=\pi^{\perp}\circ\bar\nu,$$
where $\bar\nu\!\in\!\Ga(\cZ_{\T};\fO)$ is the section induced by $\nu$.
The homomorphism~\e_ref{fOsurj_e} extends to a surjective homomorphism from 
the cokernel bundles over~$\bar\cZ_{\T}$. 
In the next two paragraphs, we show that $\fO^{\perp}\!\lra\!\bar\cZ_{\T}$ 
contains a trivial $C^\infty$-subbundle unless 
$|\Ver|\!=\!1$.
Therefore, $\tC_{\T}(*)\!=\!0$ if $\cZ_{\T}$ is not of type~\ref{basiceff_item}.

Suppose first that $(\Ver,\Edg)$ contains a loop $L\!\subset\!\Ver$ 
and $\ka$ is not constant on $L$.
Then, the image of the principal components $\Si_P$ 
of any element $[\Si,u]$ of $\cZ_{\T}$ contains
at least two curves in~$X$.
Then, $\fO^{\perp}$ contains a pull-back of the bundle
$$\E^*\otimes f^*T\ov\cM_{\un\be}\oplus s^*\N_{|\un\be|}\lra
\ov\cM_{1,1}\times\ov\cM_{\un\be},$$
where $s\!:\ov\cM_{\un\be}\!\lra\!\ov\cM_{\un\be,1}$ is the bundle section
taking each curve $\cC$ to one of the nodes.
The bundle contains a trivial $C^\infty$-subbundle for dimensional reasons.

We next consider the remaining cases.
Let $P\!\in\![m]$ be the component of the curves in $\ov\cM_{\un\be}$
containing the image of the principal component 
$\Si_P$ of any element $[\Si,u]$ of $\cZ_{\T}$.
Denote by $\ov\cM_{\un\be,1}^P\!\subset\!\ov\cM_{\un\be,1}$ the component
consisting of the curves~$\cC_P$, with $\cC\!\in\!\ov\cM_{\un\be}$.
If $X$ is an ideal Calabi-Yau
 $5$-fold and $m\!>\!1$ (implying $\cC\!\neq\!\cC_P$),
the restriction of $\N_{|\un\be|}$ to $\cC_P\!\approx\!\P$ splits as either 
$\cO\!\oplus\!\cO$ or $\cO\!\oplus\!\cO(-1)$.
If $m\!=\!1$, the splitting is $\cO(-1)\!\oplus\!\cO(-1)$.
\begin{enumerate}[label={\it{Case}} \arabic*:,leftmargin=*]

\item\label{singcC_item1} $m\!>\!1$ and 
the restriction of $\N_{|\un\be|}$ to $\cC_P$ splits as $\cO\!\oplus\!\cO$.
Here, $Q|_{\ov\cM_{\un\be,1}^P}\!=\!f^*\bar{Q}$ 
for a bundle $\bar{Q}\!\lra\!\ov\cM_{\un\be}$.
Since the restriction of $f^*\bar{Q}$ to $u(\Si_P)$ is trivial,
$\fO^{\perp}$ contains the subbundle $\E^*\!\otimes\!f^*\bar{Q}$,
where $\E\!\lra\!\ov\M_1(X,\be)$ is the Hodge line bundle.
The subbundle is a pull-back of the bundle 
$$\E^*\otimes\bar{Q}\lra\ov\cM_{1,1}\times\ov\cM_{\un\be},$$
which contains a trivial $C^\infty$-subbundle
for dimensional reasons by \e_ref{Qdfn_e}.

\item\label{singcC_item2} $m\!>\!1$ and 
the restriction of $\N_{|\un\be|}$ 
to $\cC_P$ splits as $\cO\!\oplus\!\cO(-1)$. 
Here, $Q|_{\ov\cM_{\un\be,1}^P}$ contains a subbundle $f^*\bar{Q}'$ 
of co-rank $1$ for a bundle $\bar{Q}'\!\lra\!\ov\cM_{\un\be}$.
Since the restriction of $f^*\bar{Q}'$ to $u(\Si_P)$ is trivial,
$\fO^{\perp}$ contains the subbundle $\E^*\!\otimes\!f^*\bar{Q}'$,
which is a pull-back of the bundle
$$\E^*\otimes\bar{Q}'\lra\ov\cM_{1,1}\times\ov\cM_{\un\be}.$$
Thus, the subbundle admits a section~$s$ such that $s^{-1}(0)$
is contained in the union of the spaces $\ov\M_1(\cC_i,A_i)$  
taken over finitely many $\un\be$-curves~$\cC_i$.
Since the restriction of~$\N_{|\un\be|}$ to $\cC_P$ contains $\cO(-1)$, 
$\fO^{\perp}$ also contains a line subbundle is isomorphic 
to $\ev_{z_1}^*L$, where $\ev_{z_1}\!:\bar\cZ_{\T}\!\lra\!X$ is the evaluation
map sending $[\Si,u]$ to the value of $u$ at a node of~$\Si$
taken to a node of~$\cC_P$.
The restriction of this subbundle to $s^{-1}(0)$ is trivial.

\item\label{smoothcC_item} $m\!=\!1$.
Here, $\fO^{\perp}\!=\!\fO$ is a bundle of the same rank 
as the dimension of $\ov\M_1(\ov\cM_{\be,1},d)$ for some $d\!\in\!\Z^+$.
Thus, if $\cZ_{\T}$ is not the main stratum of $\ov\M_1(\ov\cM_{\be,1},d)$, 
the restriction of $\fO^{\perp}$ to $\bar\cZ_{\T}$ contains a 
trivial $C^{\i}$-subbundle.\\
\end{enumerate}

It remains to consider the case $\cZ_{\T}\!=\!\cZ_{\T_{\eff}(\be,d)}$
with $|\Ver|\!=\!1$.
Then, $\vph_{t\nu}$ is a generic section~of
$$\fO=\fO^{\perp}=\E^*\otimes f^*T\ov\cM_{\be}\oplus R^1\pi_*\ev^*\N_{\be}
\lra\ov\M_1^0(\ov\cM_{\be,1},d),$$
where $\ov\M_1^0(\ov\cM_{\be,1},d)\!\subset\!\ov\M_1(\ov\cM_{\be,1},d)$
is the closure of the space of maps with smooth domains,
$\pi$ is the structure map for the universal curve over $\ov\M_1(\ov\cM_{\be,1},d)$,
and $\ev$ is the corresponding evaluation map, see~\e_ref{univcurv_e}.
Thus, 
\BE{effcontr_e1}\begin{split}
\tC_{\T_{\eff}(\be,d)}(d\be)
&=\blr{e\big(\E^*\otimes f^*T\ov\cM_{\be}\big)e\big( R^1\pi_*\ev^*\N_{\be}\big), 
\ov\M_1^0(\ov\cM_{\be,1},d)}\\
&=\blr{c_2\big(\ov\cM_{\be}\big),\ov\cM_{\be}}
\int_{\ov\M_1^0(\P,d)}e\big(R^1\pi_*\ev^*\big(\cO(-1)\!\oplus\!\cO(-1)\big)\big)\\
&\qquad-\blr{\la\,f^*c_1\big(\ov\cM_{\be}\big)e\big( R^1\pi_*\ev^*\N_{\be}\big), 
\ov\M_1^0(\ov\cM_{\be,1},d)},
\end{split}\EE
where $\la\!=\!c_1(\E)$.
Using the Atiyah-Bott Localization Theorem of~\cite{AtB} 
as in Section~27.5 of~\cite{MirSym}, we find 
\BE{effcontr_e3}
\int_{\ov\M_1^0(\P,d)}e\big(R^1\pi_*\ev^*\big(\cO(-1)\!\oplus\!\cO(-1)\big)\big)
=\frac{d-1}{12d^2} ~.\footnotemark\EE
In the next paragraph, we will obtain 
\BE{effcontr_e4}\begin{split}
&\blr{\la\,f^*c_1\big(\ov\cM_{\be}\big)e\big( R^1\pi_*\ev^*\N\big), 
\ov\M_1^0(\ov\cM_{\be,1},d)}\\
&\qquad\qquad\qquad\qquad=\frac{d-1}{24d^2}
\int_{\ov\cM_{\be,1}}c_2(\N_{\be})\,f^*c_1(\ov\cM_{\be}).
\end{split}\EE
Along with \e_ref{effcontr_e1}, \e_ref{effcontr_e3}, and~\e_ref{ghcontr_e},
we conclude the first identity in~\e_ref{g1thm_e}.

\footnotetext{The space $\ov\M_1^0(\P,d)$ is singular, and thus 
\cite{AtB} is not generally applicable.
However, with the setup of Section~27.5 in~\cite{MirSym},
the restrictions of 
$e\big(R^1\pi_*\ev^*\big(\cO(-1)\!\oplus\!\cO(-1)\big)\big)$
to all fixed loci, with the exception of the simplest 1-edged ones, vanish.
Furthermore, $\ov\M_1^0(\P,d)$ is smooth along the 1-edged fixed loci.
Therefore, the usual Atiyah-Bott localization formula applies.
The normal bundles to the only contributing loci are the same 
as the normal bundles in the desingularization $\wt\M_1^0(\P,d)$
of $\ov\M_1^0(\P,d)$ constructed in~\cite{VaZ} and described
in Section~1.4 in~\cite{VaZ}.}

With $\M_{0,2}(\ov\cM_{\be,1},d)\!\subset\!\ov\M_{0,2}(\ov\cM_{\be,1},d)$
denoting  the locus of maps with nonsingular domains, let 
$$\M_{0,1=2}(\ov\cM_{\be,1},d)
=\big\{b\!\in\!\M_{0,2}(\ov\cM_{\be,1},d)\!: \ev_1(b)\!=\!\ev_2(b)\big\}.$$
Denote by $\ov\M_{0,1=2}^0(\ov\cM_{\be,1},d)$ the closure of 
$\M_{0,1=2}^0(\ov\cM_{\be,1},d)$ in $\ov\M_{0,2}(\ov\cM_{\be,1},d)$.
Let 
$$\ov\De_1^0(\ov\cM_{\be,1},d)\subset\ov\M_1^0(\ov\cM_{\be,1},d)$$
be the subspace consisting of the stable maps $[\Si,u]$
such that the principal component $\Si_P$ of $\Si$ is singular. 
There is a natural node-identifying immersion
$$\io\!:\ov\M_{0,1=2}^0(\ov\cM_{\be,1},d)\big/S_2
\lra \ov\De_1^0(\ov\cM_{\be,1},d),$$
which is an embedding outside of a divisor.
Note  
$$\io^*e\big( R^1\pi_*\ev^*\N_{\be}\big)
=c_2(\N_{\be})\,e\big(R^1\pi_*\ev^*\N_{\be}\big).$$
If $\De_1\!\subset\!\ov\cM_{1,1}$ is the locus of the nodal elliptic
curve,
$$\la=\frac{1}{12}\De_1\in H^2\big(\ov\cM_{1,1}\big).$$
Therefore,
\BE{effcontr_e5}\begin{split}
&\blr{\la\,f^*c_1\big(\ov\cM_{\be}\big)e\big( R^1\pi_*\ev^*\N_{\be}\big), 
\ov\M_1^0(\ov\cM_{\be,1},d)}\\
&\qquad =\frac{1}{24}
\blr{f^*c_1\big(\ov\cM_{\be}\big)\,c_2(\N_{\be})
e\big(R^1\pi_*\ev^*\N_{\be}\big),\ov\M_{0,1=2}^0(\ov\cM_{\be,1},d)}\\
&\qquad=\frac{1}{24}
\blr{c_2(\N_{\be})f^*c_1\big(\ov\cM_{\be}\big),\ov\cM_{\be,1}}\\
&\qquad\qquad\qquad\qquad\times
\int_{\ov\M_{0,1=2}^0(\P_p,d)}
e\big(R^1\pi_*\ev^*\big(\cO(-1)\!\oplus\!\cO(-1)\big)\big).
\end{split}\EE
Using the localization as in Section~27.5 of~\cite{MirSym} once again, 
we find 
$$\int_{\ov\M_{0,1=2}^0(\P_p,d)}
e\big(R^1\pi_*\ev^*\big(\cO(-1)\!\oplus\!\cO(-1)\big)\big)
=\frac{d-1}{d^2} ~.\footnotemark$$
Along with \e_ref{effcontr_e5}, this identity implies \e_ref{effcontr_e4}.

\footnotetext{Similarly to the situation discussed for \e_ref{effcontr_e3}, 
$\ov\M_{0,1=2}^0(\P_p,d)$ has singularities but is nonsingular along
the only fixed locus to which 
$e\big(R^1\pi_*\ev^*\big(\cO(-1)\!\oplus\!\cO(-1)\big)\big)$
restricts non-trivially.}

\section{Local $\bP^2$}

\subsection{Gromov-Witten invariants}
\label{lll1}

\noindent
We consider here the local Calabi-Yau $5$-fold given
by the total space
\BE{localP2dfn_e}
X= \cO(-1) \oplus \cO(-1) \oplus \cO(-1)\lra \bP^2 \,.\EE
There are only two primary Gromov-Witten invariants in each degree $d$:
$$N_{0,d} = N_{0,d}(H^2,H^2) \ \ \text{ \ and\  } \  \
N_{1,d},$$ where $H$ is the hyperplane class in $H^2(X,\mathbb{Z})\!=
\!H^2(\bP^2,\mathbb{Z})$.
We compute both Gromov-Witten invariants by localization{\footnote{In the
genus 0 case, the moduli space $\ov\M_{0,2}(\bP^2,d)$
is a nonsingular
 stack and
the usual Atiyah-Bott localization formula applies. In the
genus 1 case, the virtual localization formula of \cite{GrP} is used.}} and then state 
a conjectural formula found by Martin for the integer counts $n_{1,d}$.

\begin{lmm}\label{localP2_lmm} For $d\!\in\!\Z^+$,
$$N_{0,d}= \frac{(-1)^{d-1}}{d} \qquad\hbox{and}\qquad
N_{1,d} = \frac{(-1)^d}{8d} \,.$$
\end{lmm}

\noindent {\em Proof.}
Let $(a,b,c)$ be the weights of the torus action on the
vector space $\C^3$.
The weights of the torus action on~$T\bP^2$ at the fixed points are then
\begin{equation*}\begin{split}
&P_1\!=\![1,0,0]:\qquad      b\!-\!a, c\!-\!a,\\
&P_2\!=\![0,1,0]:\qquad      a\!-\!b, c\!-\!b,\\
&P_3\!=\![0,0,1]: \qquad     a\!-\!c, b\!-\!c.
\end{split}\end{equation*}
We choose linearizations on the 3 bundles $O(-1)$ with the
following weights at the fixed points:
\begin{equation*}
\begin{array}{cccc}
&\cO(-1)&    \cO(-1)&   \cO(-1)\\
P_1:\qquad& 0 &       a\!-\!b&     a\!-\!c\\
P_2:\qquad&        b\!-\!a&      0&       b\!-\!c\\
P_3:\qquad&        c\!-\!a&      c\!-\!b&     0\\
\end{array}
\end{equation*}

In order to compute the numbers $N_{0,d}$,
we choose the points $P_1$ and $P_2$ for the insertions and integrate over
$$\ov\M=\big\{b\!\in\!\ov\M_{0,2}(\bP^2,d)\!: 
\ev_1(b)\!=\!P_1,~\ev_2(b)\!=\!P_2\big\}.$$
By the choice of the weights and the points, there is 
a unique fixed locus with non-zero contribution,
see Section 27.5 in \cite{MirSym} for a similar situation.
The locus consists of the $d$-fold cover $u$ of the line 
$$\bP^1_{12}\!=\!\ov{P_1P_2}$$ 
branched over only $P_1$
and~$P_2$ and with the marked points $1$ and $2$ mapped to
$P_1$ and~$P_2$, respectively.
The weights of the fibers of the relevant bundles at the fixed locus are given~by
\begin{equation*}\begin{array}{ll}
H^1(u^*\cO(-1)):\qquad&      \frac{(-1)^{d-1}(d-1)!}{d^{d-1}}  (a-b)^{d-1}\\
\noalign{\medskip}
H^1(u^*\cO(-1)):\qquad&       \frac{(-1)^{d-1}(d-1)!}{d^{d-1}}  (b-a)^{d-1}\\
\noalign{\medskip}
H^1(u^*\cO(-1)):\qquad&      (-1)^{d-1}\prod_{r=1}^{d-1}\Big(c-\frac{(d-r)a+rb}{d}\Big)\\
\noalign{\medskip}
T\ov\M:&               \frac{(-1)^{d-1}(d-1)!^2}{d^{2(d-1)}}(a\!-\!b)^{2(d-1)}
\prod_{r=1}^{d-1}\Big(c-\frac{(d-r)a+rb}{d}\Big)\,,\\
\end{array}\end{equation*}
see Section 27.2 in \cite{MirSym}.
The number $N_{0,d}$ is the ratio of the product of the first three expressions
and the last expression, divided by $d$ for the stack automorphism factor.

We next compute the number $N_{1,d}$.
There are now $6$ fixed loci with nonzero contribution: 
the three $d$-fold Galois covers of the three lines together 
with a choice of vertex for the contracted elliptic component.
By symmetry, the contribution of the $d$-fold cover of $\bP^1_{12}$
with the contracted component at $P_1$ determines the
other cases.
The weights of the fibers of the relevant bundles at 
the  $d$-fold cover of $\bP^1_{12}$  are given by
\begin{equation*}\begin{array}{ll}
H^1(u^*\cO(-1)):\qquad&  \frac{(-1)^{d-1}(d-1)!}{d^{d-1}} (a-b)^{d-1}(-\la)\\
\noalign{\medskip}
H^1(u^*\cO(-1)):\qquad&   \frac{(-1)^{d-1}(d-1)!}{d^{d-1}}(b-a)^{d-1}(a-b-\la)\\
\noalign{\medskip}
H^1(u^*\cO(-1)):\qquad&      (-1)^{d-1}\prod_{r=1}^{d-1}\Big(c-\frac{(d-r)a+rb}{d}\Big)
\,(a-c-\la)\\
\noalign{\medskip}
Obs(\bP^2):&  (b-a-\la)(c-a-\la)\\
\noalign{\medskip}
T\ov\M_1(\bP^2,d):& \frac{(-1)^d d!^2}{d^{2d-1}}(a\!-\!b)^{2d-1}
\prod_{r=0}^{r=d}\Big(c-\frac{(d-r)a+rb}{d}\Big)
\,\Big(\frac{b-a}{d}-\psi\Big),\\
\end{array}\end{equation*}
where $\la$ is the first chern class of the Hodge line bundle 
$\E\!\lra\!\ov\cM_{1,1}$.
The contribution of the locus to $N_{1,d}$ is 
the ratio of the product of the first four expressions
and the last expression, divided by the stack factor $d$, 
and integrated over~$\ov\cM_{1,1}$,
$$\Cont(a,b)=\frac{(-1)^d}{24d}\,\frac{c-a}{c-b}\,.$$
Symmetrizing over $a$, $b$, and $c$, we obtain $N_{1,d}$.\QED \\
\vspace{+10pt}

By Lemma \ref{localP2_lmm} and the $n\!=\!2$ case of \e_ref{g0BPS_e},
the genus $0$ counts for $X$ are given by
$$n_{0,d}=n_{0,d}(H^2,H^2)
=\begin{cases} 1,&\hbox{if}~d\!=\!1;\\
-1,&\hbox{if}~d\!=\!2;\\
0,&\hbox{if}~d\!\ge\!3.
\end{cases}$$
Using the algorithm of Section \ref{dim5CY_subs}, 
we have computed the genus $1$ count $n_{1,d}$ for $X$ for $d\!\le\!200$.
All are integers.

\subsection{Martin's conjecture}
\label{lll2}
Recall the definition of the M\"obius $\mu$-function,
$$\mu\!:\Z^+\lra\{0,\pm1\}, \qquad
\mu(d)=\begin{cases}
(-1)^r,&\hbox{if}~d~\hbox{is the product of $r$ distinct primes},\\ 
0,&\hbox{otherwise}.\end{cases}$$
Define a sign function $S(d)$ and an absolute value function $V(d)$ as
follows:
$$S(d)=\begin{cases}
\mu(d),&\hbox{if}~d\not\cong 4~(\mod 8),\\
\mu(d/4),&\hbox{if}~d\cong 4~(\mod 8),
\end{cases}\qquad
V(d)=\frac{k^2\!-\!1}{8}\times
\begin{cases}
\frac{k^2-1}{8},&\hbox{if}~d\!=\!k,\,2\not|k,\\
\frac{17k^2+7}{8},&\hbox{if}~d\!=\!2k,\,2\not|k,\\
2k^2\!+\!1,&\hbox{if}~d\!=\!4k,\,2\not|k.
\end{cases}$$

\begin{cnj}[G.~Martin]
\label{localP2_cnj}
For every $d\!\in\!\Z^+$, the genus~$1$ degree~$d$ count 
for the local Calabi-Yau $5$-fold $\bP^2$ is given by
\BE{localP2n1_e} n_{1,d}=S(d)V(d).\EE
\end{cnj}

If $8|d$, then $S(d)$ vanishes and a definition of 
$V(d)$ is not required for \eqref{localP2n1_e}.
As our method for computing the numbers $n_{1,d}$ from $n_{0,d}$ and 
$N_{1,d}$ is completely explicit and the starting data is fairly simple,
a verification of Conjecture \ref{localP2_cnj} by elementary
identities may be possible. 
Unfortunately, the algorithm involves a 
significant number of simultaneous recursions.\footnote{Explicit forms
of these recursions can be found in the appendix to this paper available
from the authors' websites.}

Geometric consequences are easily obtained from the conjecture.
For example, since $n_{1,d}$ is predicted to vanish whenever
$8|d$, we expect Calabi-Yau 5-folds obtained from
suitably generic deformations of the local $\mathbb{P}^2$ geometry
to contain {\em no embedded elliptic curves of degrees divisible
by~8}. Is there a simple symplectic reason for this?

\vspace{+10 pt}
\noindent
Department of Mathematics\\
Princeton University\\
Princeton, NJ 08544, USA\\
rahulp@math.princeton.edu

\vspace{+10 pt}
\noindent
Department of Mathematics \\
SUNY Stony Brook\\
Stony Brook, NY 11794, USA\\
azinger@math.sunysb.edu \\

\end{document}